\pdfoutput=1
\documentclass[pdflatex,sn-mathphys,Numbered]{sn-jnl}

\usepackage{graphicx}%
\usepackage{multirow}%
\usepackage{amsmath,amssymb,amsfonts}%
\usepackage{amsthm}%
\usepackage{mathrsfs}%
\usepackage{mathtools}
\usepackage{enumitem}
\usepackage{xcolor}%
\usepackage{textcomp}%
\usepackage{manyfoot}%
\usepackage{booktabs}%
\usepackage{algorithm}%
\usepackage{algpseudocodex}%
\usepackage{listings}%
\usepackage[utf8]{inputenc}
\usepackage{color}
\usepackage{nicefrac}
\usepackage{comment}
\usepackage{multirow,multicol,longtable,makecell,booktabs}
\usepackage[numbers]{natbib}

\usepackage[font=small, labelfont=bf, justification=RaggedRight, font=small]{caption}
\usepackage{xfrac}	
\usepackage[capitalise,noabbrev]{cleveref}	
\crefformat{equation}{#2(#1)#3}
\crefrangeformat{equation}{(#3#1#4) to (#5#2#6)}
\crefmultiformat{equation}{#2(#1)#3}{ and~#2(#1)#3}{, #2(#1)#3}{, and~#2(#1)#3}

\usepackage {tikz}
\usepackage{tikzsymbols}

\def \NP {$\mathcal{NP}$}

\allowdisplaybreaks

\newcommand{\dsum}{\displaystyle \sum}
\newcommand{\sumin}[1]{\sum_{#1}}
\newcommand{\dsumin}[1]{\dsum_{#1}}
\newcommand{\xVar}[1][i,j]{x_{#1}^{e,s}}
\newcommand{\yVar}[1][i,j]{y_{#1}^{k,s}}
\newcommand{\zVar}[1][i,j]{z_{#1}^{e,s}}
\newcommand{\uVar}[1][s]{u^{#1}}

\newcommand{\fVar}[1][i]{f_{#1}^{e,s}}
\newcommand{\bVar}[1][i]{b_{#1}^{e,s}}
\newcommand{\wVar}[1][i]{w_{#1}^{e,s}}
\newcommand{\BVar}[1][i]{B_{#1}}

\newcommand{\problema}{\emph{IVPRP}}
\newcommand{\modelo} {\emph{IVPRM}}
\newcommand{\Rmodelo}{\emph{R-IVPRM}}

\newcommand{\orcidlink}[1]{%
\textbf{}
}

\definecolor{siamc}{HTML}{00A499}

\usepackage{geometry}
\theoremstyle{thmstyleone}%
\theoremstyle{thmstyletwo}%
\newtheorem{remark}{Remark}%

\theoremstyle{thmstylethree}%

\usepackage{xkcdcolors}
\usepackage{hyperref}
\hypersetup{
    colorlinks = true,
    urlcolor = siamc,
    linkcolor = xkcdMagenta,
    citecolor = xkcdLeafGreen
}

\begin{document}

\title[]{The integrated vehicle and pollster routing problem
\newline{\color{white}\fontsize{0.1pt}{0.1pt}\selectfont \emph{$\star$ 10-year anniversary $\star$}}
\vspace{-1.5em}
}


\author[1]{\fnm{Sandra} \sur{Gutiérrez} \orcidlink{0000-0002-0409-4601}}\email{sandra.gutierrez@epn.edu.ec}
\equalcont{These authors contributed equally to this work.}

\author*[2,3,4]{\fnm{Andrés} \sur{Miniguano-Trujillo} \orcidlink{0000-0002-0877-628X}}\email{andres.miniguano-trujillo@kcl.ac.uk}
\equalcont{These authors contributed equally to this work.}

\author[1,3]{\fnm{Diego} \sur{Recalde} \orcidlink{0000-0002-6981-2272}}\email{diego.recalde@epn.edu.ec}
\equalcont{These authors contributed equally to this work.}

\author[3]{\fnm{Luis M} \sur{Torres} \orcidlink{0000-0002-0778-195X}}\email{luis.torres@epn.edu.ec}
\equalcont{These authors contributed equally to this work.}

\author[1]{\fnm{Ramiro} \sur{Torres} \orcidlink{0000-0003-2057-1719}}\email{ramiro.torres@epn.edu.ec}
\equalcont{These authors contributed equally to this work.}

\affil[1]{\orgdiv{Department of Mathematics}, \orgname{Escuela Polit\'ecnica Nacional}, \orgaddress{\street{Ladr\'on de Guevara}, \city{Quito}, \postcode{170525}, \state{Pichincha}, \country{Ecuador}}}

\affil[2]{\orgdiv{Maxwell Institute for Mathematical Sciences}, \orgaddress{\street{Bayes Centre}, \city{Edinburgh}, \country{United Kingdom}}}

\affil[3]{\orgdiv{Research Center for Mathematical Modeling -- ModeMat}, \orgname{Escuela Polit\'ecnica Nacional}, \orgaddress{\street{Ladr\'on de Guevara}, \city{Quito}, \postcode{170525}, \state{Pichincha}, \country{Ecuador}}}

\affil[4]{\orgdiv{King's College London}, \orgaddress{\street{Centre for Neuroimaging Sciences}, \city{London}, \country{United Kingdom}}}



\abstract{The National Statistics Bureau of Ecuador conducts monthly price surveys to monitor the evolution of the Consumer Price Index, a metric that measures the consumer prices of essential commodities. These surveys are administered across a designated set of stores, with a fleet of vehicles transporting pollsters from the bureau headquarters to the chosen locations. Additionally, pollsters move between stores using walking paths or vehicles to reduce travel time. This paper introduces the Integrated Vehicle and Pollster Routing Problem and presents an integer programming model to effectively schedule pollster visits to selected stores while optimizing the routing of the vehicle fleet. A three-phase algorithm is proposed, and computational experience with real-world instances is presented.}

\keywords{vehicle routing, crew scheduling, integer programming, primal heuristics}


\maketitle

\section{Introduction}\label{sec:Introduction}

This paper addresses an original problem that incorporates several elements from the Simultaneous Vehicle Routing and Crew Scheduling Problem (SVRCSP). The problem stems from data collection procedures implemented by the National Statistics Bureau of Ecuador (INEC) for computing the national Consumer Price Index (CPI). A team of pollsters and a fleet of vehicles traverse walking paths  and vehicular routes to visit stores scattered throughout the city of Guayaquil and its outskirts over several days. On each workday, a pollster visits each store to gather prices of specific goods necessary to calculate the CPI. Additionally, pollsters have the option to walk or be transported by fixed-capacity vehicles during their shifts, which consist of morning and afternoon periods separated by a mandatory break. Each vehicular route must commence at INEC headquarters, transporting pollsters and vehicles to and from stores, allowing for vehicle returns or waiting at each location. At the end of each workday, vehicles and pollsters must return to the headquarters. The primary objective of this problem is to minimize the total fixed daily costs, which include expenses related to depot maintenance, pollster hire, and vehicle hire.

The problem addressed in this work is closely related to recent studies on the SVRCSP, which combines vehicle routing and crew scheduling. In this class of problems, designing vehicle routes and assigning and moving crews are interconnected decisions that must be made jointly. Such integration is relevant in real-world applications where crews operate at multiple locations, may transfer between vehicles, and need to meet strict synchronization requirements. These interactions considerably increase the computational complexity of the problem and set it apart from treating vehicle routing and crew scheduling as separate planning tasks. From a methodological standpoint, the literature on the SVRCSP can be broadly divided into two streams: heuristic approaches, which are the most common in practice due to the size and complexity of real instances, and exact approaches, which provide stronger guarantees but are typically limited to smaller or more structured settings.

In the SVRCSP literature, the most common approach relies on heuristic methods based on a two-stage sequential decomposition. Typically, these methods construct vehicle routes in the first stage, often without considering driver constraints, and assign drivers to route segments in the second stage while ensuring compliance with labor regulations. In this vein, \citet{Fikar2015} examine a project motivated by the Austrian Red Cross, in which the authors address the daily planning of home healthcare providers operating multiple vehicles to transport nurses to clients and retrieve them after service. Their solution procedure is a matheuristic consisting of two phases: identifying potential walking routes and optimizing the transport system. \citet{Molenbruch2021} investigate an integrated mobility system in which a user's trip may combine dial-a-ride services with regular public transportation. The authors develop a routing algorithm and an integrated scheduling procedure to synchronize flexible vehicle routes with public transport timetables, thereby supporting both the design and operational implementation of the system. \citet{Drexl2013} study a simultaneous vehicle and crew routing and scheduling problem arising in long-distance road transport in Europe, where pickup-and-delivery requests must be fulfilled over a multi-period planning horizon using a heterogeneous fleet of trucks and drivers. The authors devise a heuristic solution method based on a two-stage decomposition while accounting for EU social legislation governing drivers. \citet{Lucci2021} present a simultaneous vehicle routing and crew scheduling problem in which long-distance pickup-and-delivery requests must be fulfilled over a multi-day planning horizon, subject to constraints such as multiple time windows and hours-of-service regulations. A two-stage sequential approach is adopted: a set of truck routes is generated in the first stage, and a set of driver routes consistent with the truck routes is determined in the second stage. The authors develop an algorithm based on the GRASP \(\times\) ILS metaheuristic, embedded with a repair procedure to facilitate the construction of feasible initial solutions. \citet{Mendes2020} describe a decision support system developed for an Italian pharmaceutical distribution company to address a multi-trip vehicle routing problem (VRP) with additional truck and driver scheduling constraints. The authors propose a two-phase algorithm that applies an Iterated Local Search metaheuristic to design vehicle routes and a mixed-integer programming (MIP) model to assign trucks and drivers to those routes.

Exact methods, despite their computational cost, yield rigorous bounds and solutions that shed light on the combinatorial structure of a problem of interest. Thus, they are highly relevant for evaluating the quality of heuristic strategies. As a consequence, some studies have explored this trend. \citet{Goel2017} study the combined vehicle routing and truck driver scheduling problem by considering working-hour constraints. The authors proposed a branch-and-price algorithm, which employs a dynamic programming-based labeling algorithm for the generation of routes complying with the regulations. \citet{DomnguezMartn2018} introduce a multi-depot SVRCSP in which vehicles depart from one depot and arrive at another, while drivers both start and end their routes at the same depot. Driver routes are circular and subject to maximum-duration constraints, whereas vehicle routes are non-circular. Motivated by an air-traffic application, the problem is formulated as an integer program (IP) strengthened with valid inequalities and solved using an exact algorithm. More flexible driver exchanges are studied by \citet{Ammann2023}, where drivers may transfer between buses at intermediate stops, introducing additional synchronization constraints. The authors propose a mathematical formulation based on a time-expanded multidigraph, derive bounds on the minimum number of drivers required, and develop a destructive-bound-enhanced matheuristic that converges to provably optimal solutions. Their approach is validated through a real-world case study involving a bus company. \citet{Lam2020} present a setting in which crews can interchange vehicles and travel as passengers, creating spatial and temporal interdependencies between vehicle and crew routes. The authors propose three solution approaches: a constraint programming model that overlays crew-routing constraints onto a standard vehicle routing problem, an MIP formulation solved directly with an IP solver, and a two-stage method. Their results show that single-stage approaches outperform sequential methods and that decoupling crews from vehicles improves solution quality.
\citet{Lucci2025} investigate the SVRCSP in long-haul transportation with time windows over a multi-day planning horizon. In this context, trucks and drivers are not permanently paired: drivers may be exchanged at designated locations and can also travel as passengers or use external taxi services at an additional cost. The objective is to minimize truck operating costs, taxi expenses, and penalties for late deliveries. The authors propose three compact MIP formulations and introduce several families of valid inequalities employed as cutting planes within a branch-and-cut algorithm.

While several of the cited works consider multi-period settings with simultaneous vehicle routing and crew scheduling, and some allow crew members to travel independently of vehicles, the problem studied here differs in several fundamental aspects that prevent direct application of existing approaches. First, \textbf{pollsters may move independently} on foot along a designated pedestrian network with fixed travel times, which directly constrain path feasibility. A pollster may be dropped off at one location by a vehicle and later picked up by a different vehicle. This \textbf{transshipment structure} introduces an extra layer of synchronization and requires a formulation capable of preventing inconsistent movements, disconnected walking segments, or cyclic pollster paths. Therefore, the need to coordinate pedestrian movements with vehicle movements is therefore one of the main modeling challenges. Second, each active pollster must take a \textbf{mandatory break within a specified time window} during each working day. This requirement is closely linked to the routing structure because the break must occur at a visited store during the pollster’s shift, not just during idle time. As a result, break feasibility must be coordinated jointly with service, walking, and vehicle movements. Finally, the optimization goal is to \textbf{minimize fixed daily operation costs}, namely the number of active days, vehicles, and pollsters required to complete the data-collection plan. This objective reflects the economics of a public-sector data-collection application and offers a different optimization perspective from those commonly adopted in related routing problems, where the focus is on variable travel costs, delay penalties, or external transport costs. Here, the focus is on resource allocation and feasibility over a multiple-day planning horizon. Overall, these features make the problem structurally different from those existing in the literature, and their combination requires a targeted modeling and solution framework, rather than a straightforward adaptation of prior approaches.

\subsection{Contribution and structure}\label{sec:Contribs}

This paper makes four contributions to the literature on synchronized transportation and routing problems.

\begin{enumerate}
	\item \textbf{Problem contribution.}
	We formally introduce the Integrated Vehicle and Pollster Routing Problem (IVPRP), a new problem class inspired by a real operational requirement in public-sector data collection. The IVPRP incorporates elements of the Simultaneous Vehicle Routing and Crew Scheduling Problem (SVRCSP) while differing from this framework in two essential ways. First, a pollster dropped off at a store need not be collected by the same vehicle, allowing transshipments through interactions between multiple vehicles. Second, pollsters may move between stores via walking paths, meaning service completion may involve both vehicular transport and pedestrian movement. These two features make synchronization central to the problem: vehicle routes, walking paths, service times, and daily shift limits require coordination within a single working-day period to construct a multi-day plan.

	\item \textbf{Modeling contribution.}
	We propose a MIP formulation that jointly captures the key operational features of the IVPRP: a two-speed network (walking and vehicular travel), a daily shift structure, mandatory breaks, and vehicle-as-transporter semantics. Specifically, vehicles are not modeled as service agents but as mobile resources transporting pollsters between service locations. Unlike related mixed-agent routing models where part of the movement structure is fixed in advance, our formulation treats both pedestrian and vehicular movements as optimization decisions. This results in fully synchronized routing plans and provides a modeling framework that can be extended, for example, to more complex labor-regulation settings with multiple breaks.

	\item \textbf{Algorithmic contribution.}	
	A three-phase decomposition heuristic is introduced as an initial practical solution method for this problem. Its design rationale addresses the computational challenges posed by the scale of real-world instances and is supported by a set of Integer Programming tools that together form a comprehensive solution framework. The rationale involves dividing the planning task into balanced clusters of stores that represent partial shift-coverage areas, producing bounds and solving a reduced routing problem for each cluster, and then combining feasible partial plans into a multi-day solution.

	\item \textbf{Empirical contribution.}
	We report computational results on simulated instances and a real-world case study with more than 800 stores. The experiments show that the proposed approach can solve instances that are well beyond the size that is tractable for the full integrated model, and that it can generate substantial savings in operational costs, achieving reductions of up to 50\% in a real application. 
	
\end{enumerate}

The organization of this paper is as follows: \cref{section_modelo} introduces an IP-based formulation for the \emph{Integrated Vehicle and Pollster Routing Problem}, along with symmetry-breaking inequalities and general bounds. Then, in
\cref{section_solution_approach}, the three-phase heuristic is presented alongside a set of specialized bounds. Subsequently, a real-world instance and computational results are reported in \cref{section_comp_exp}. Finally, concluding remarks are provided in \cref{section_conclusions}.

\section{An integrated model for routing pollsters and vehicles}
\label{section_modelo}

\subsection{Notation}\label{ssec:notation}

\begin{figure}[h]
	\centering
	\includegraphics[scale=0.75]{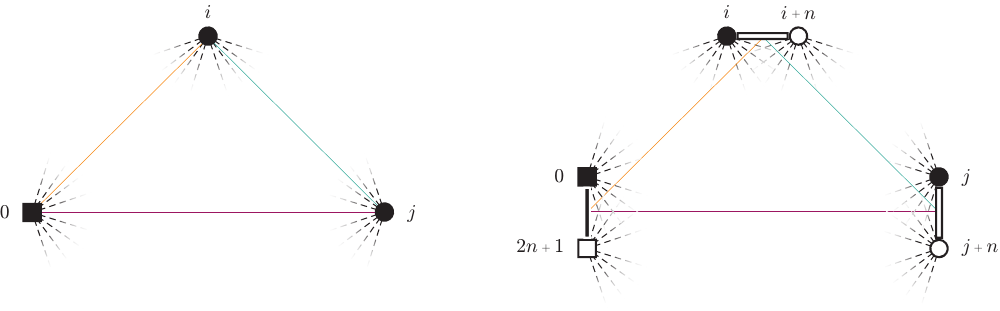}
	\caption{%
		Node duplication: The depot, labeled \(0\), and two stores, labeled \(i\) and \(j\), are associated with their duplicates, labeled \(2n+1\), \(i+n\), and \(j+n\), respectively.
	}
	\label{fig_node_duplication1}
\end{figure}

Let \(n \in \mathbb{N}\) denote the number of stores, and let \(S\) denote the set of days within the planning horizon. In this problem, only one depot is considered due to technical requirements. It is assumed that pollsters must depart from and return to the depot in a vehicle on a daily basis. Define \(G \coloneqq (V, A)\) as a multi-graph with a node set given by \(V \coloneqq \{0, \ldots, 2n+1\}\) and an arc set \(A \subseteq V \times V\). Each store and the depot are associated with two nodes in \(V\), as illustrated in \cref{fig_node_duplication1}. Specifically, the depot is represented by node \(0\), and its duplicate is node \(2n+1\). Similarly, the subsets \(C_{-} \coloneqq \{1, \dots, n\}\) and \(C_{+} \coloneqq \{n+1, \dots, 2n\}\) correspond to the stores and their duplicates, respectively. Node duplication enables the modeling of vehicle arrival, service, waiting, and departure times separately. Furthermore, define \(C \coloneqq C_{-} \cup C_{+}\), which will collectively be referred to as the set of stores.

\begin{figure}[htb!]
	\centering
	\includegraphics[scale=1.0]{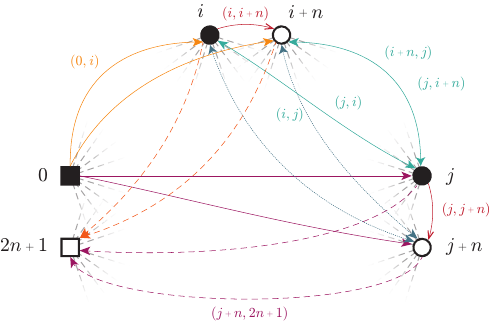}
	\caption{
	Multigraph construction for the {\problema}. Focusing only on the depot and two nodes \(i, j \in C_{-}\): there are only outgoing arcs from the depot \(0\) and only incoming arcs to its duplicate \(2n+1\). A single service arc is defined for each pair of nodes \(i\) and \(j\). Additionally, two walking arcs connect the set \(\{i, j, i+n, j+n\}\), and 18 vehicle arcs connect the set \(\{0, i, j, i+n, j+n, 2n+1\}\).
	}
	\label{fig_multigrafo2}
\end{figure}

The set \(A \coloneqq A_S \cup A_W \cup A_V\) consists of three distinct types of arcs. The first group, \(A_S \coloneqq \big\{(i, i+n) : i \in C_{-}\big\}\), comprises arcs that establish connections between each store and its corresponding duplicate, designated as service arcs. The second group of arcs, \(A_W \coloneqq \big\{(i, j) : i \in C_{+}, j \in C_{-}, i \neq j+n\big\}\), called walking arcs, links every pair of stores. Finally, the vehicle transportation arc set, \(A_V \coloneqq \big\{(i, j) : i, j \in C, i \neq j, i \neq j+n\big\} \cup \big\{(0, i) : i \in C \big\} \cup \big\{(i, 2n+1) : i \in C \big\}\), represents potential connections between nodes or depots via vehicles. Furthermore, for each arc \(a \in A_S\), the parameter \(t_a > 0\) denotes the duration of time spent by a pollster collecting data at a store. Similarly, for each arc \(a \in A_W\), \(t_a > 0\) represents the walking travel time between stores. At last, for each arc \(a \in A_V\), \(\tau_a \geq 0\) is the vehicular travel time between two nodes (\(\tau_a = 0\) only in the case of service arcs; i.e., \(a \in A_S\)). A visual representation of the resulting multigraph is provided in \cref{fig_multigrafo2}.

Visits are carried out by a set \(E\) of pollsters, transported by a homogeneous fleet of vehicles \(K\) with a positive integer capacity \(Q\); in this context, the vehicle capacity means the number of pollsters that a vehicle can carry. Each operational day incurs a fixed cost of \(\kappa_0 \geq 0\). Likewise, there are fixed hiring costs for vehicles and pollsters, \(\kappa_1 \geq 0\) and \(\kappa_2 \geq 0\) respectively, and they operate within a daily time frame \([0, B_{\max}]\). Due to labor regulations, pollsters must take a break of \(P \geq 0\) minutes starting within the time interval \([T_0, T_1] \subset [0, B_{\max} - P]\), regardless of the total time allocated to data collection. A single pause is considered throughout this work, attending to the operational necessities of the real-world application; however, an extension to multiple breaks is possible by introducing a vector of break lengths accompanied by a set of time intervals for each pollster to make a pause. A \emph{vehicular route} is a simple path from node \(0\) to \(2n+1\), exclusively using arcs in \(A_V\). A \emph{walking path} is a simple path from \(i \in C_-\) to \(j \in C_+\) consisting of an alternating sequence of arcs in the sets \(A_S\) and \(A_W\). A  \emph{feasible service path for a pollster} is defined as a simple path from \(0\) to \(2n+1\) that connects walking paths with segments of vehicular routes. Consequently, at every node \(i\), a pollster dedicates \(t_{i,i+n}\) units of time to complete a task, adhering to the maximum daily duty length. Observe that two necessary conditions for feasibility are that \( P \leq B_{\max} - \big( \min_{(0,j) \in A_V} \tau_{(0,j)} + \min_{(j,n+1) \in A_V} \tau_{(j,n+1)} \big)  \) and \(T_1 \geq \min_{(0,j) \in A_V} (\tau_{(0,j)} + t_{j,j+n} ) \); i.e., it should be possible to take a pause within the transportation times of a working day and enough time should be allocated for it to start.
Throughout this work, it is assumed that a feasible service path always commences and finishes with a vehicular segment. This constraint is particularly motivated by the real-world application and is reasonable for servicing stores located far from the depot.
Finally, a \emph{feasible daily plan} comprises distinct feasible service paths synchronized with the durations of each walking path and vehicular routes on a single day. Any sequence of feasible daily plans forms a \emph{multi-day plan}.

\begin{quote}\normalsize
\textbf{\emph{Problem:}}
The Integrated Vehicle and Pollster Routing Problem (\problema) consists of finding an multi-day plan to visit every node \(i \in C\) precisely once at minimum cost. This plan must hire at most \(|K|\) vehicles with capacity \(Q\) and \(|E|\) pollsters within a time frame not exceeding \(|S|\) days. 
\end{quote}

As an extension of the Traveling Salesman Problem (TSP), the {\problema} is thereby {\NP}--hard. The \(\mathcal{NP}\)--hardness condition for the Integrated Vehicle and Pollster Routing Problem renders it unlikely to find an optimal solution in polynomial time. In view of this complexity, the emphasis of the paper is placed on modeling and on the design of a practical decomposition-based solution approach. Consequently, \cref{section_solution_approach} is dedicated to identifying heuristic methods that can yield good-quality solutions within a reasonable amount of time.

\subsubsection*{An illustrative example}

\begin{figure}[h]
  \centering
  \includegraphics[scale=1.0]{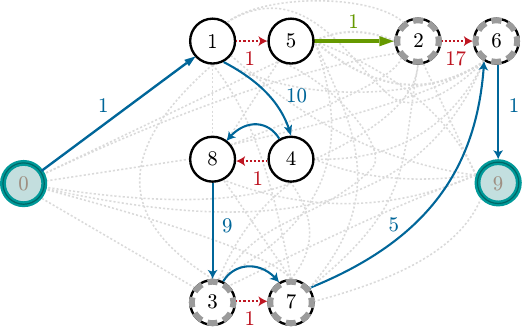}
\caption{Multigraph and optimal solution for an instance of {\problema}. }
	\label{fig_instancia_inventada}
\end{figure}

A small instance is presented in \cref{fig_instancia_inventada}, employing the following parameters: \(n = 4\) stores, \(|E| = 2\) pollster, \(|K| = 2\) vehicles with a capacity of \(Q = 2\) passengers, and a limit of \(|S| = 2\) days. The break length is fixed at one minute within a time frame of \([T_0, T_1] = [20, 30]\), and the time horizon is set to 30 minutes. The fixed and variable costs are \(\kappa_0 = \$300\), \(\kappa_1 = \$100\), and \(\kappa_2 = \$80\), respectively.

The vehicular and walking travel times are shown in \cref{table:01}. The diagonal of the pollster travel time matrix represents the service times. It is important to note that this table excludes information related to nodes in the set \(\{5, 6, 7, 8, 9\}\), which are the duplicated nodes introduced during the multigraph construction, see \cref{ssec:notation}.

\begin{table}[h]
\caption{Vehicle and walking traveling times within each pair of stores.}
\fontsize{8}{8.5}\selectfont
\setlength{\tabcolsep}{10pt}
\def\arraystretch{1.5}
\begin{tabular}{l c  c  c c c  c  c  c  c  c  c  c c}
\toprule
  & & \multicolumn{5}{c}{\bfseries Vehicle (\(\tau\))}  & \multicolumn{4}{c}{\bfseries Pollster (\(t\))} 
  \\
\cmidrule(lr){3-7} \cmidrule(lr){8-11}
& & $0$ & $1$ & $2$ & $3$ & $4$ &		 $1$ & $2$ & $3$ & $4$   \\
\hline
\textbf{Store} &$0$\, 	& -- 	  & $1$	& $4$ & $6$ & $2$				& -- & -- & -- & -- 
\\
			&$1$\, 	& $2$ & -- 	& $1$   & $2$ & $10$ 			& $1$ & $2$ & $3$ & $14$ 
			\\
			&$2$\, 	& $1$ & $1$   	& --   & $4$   & $15$ 				& $2$ & $17$ & $8$ & $20$ 
			\\
			&$3$\, 	& $7$ & $2$ 	& $5$ & -- & $9$ 				& $3$ & $8$ & $1$ & $15$ 
			\\
			&$4$\, 	& $1$ & $15$   	& $11$   & $9$   & -- 				& $15$ & $20$ & $15$ & $1$ 
			\\
\bottomrule
\end{tabular}
\label{table:01}
\end{table}
As depicted in \cref{fig_instancia_inventada}, the optimal solution for this small instance consists of a single vehicular route: \(0 \to 1 \to 4 \to 8 \to 3 \to 7 \to 6 \to 9\) (continuous blue lines), and three walking paths: \(1 \to 5 \to 2 \to 6\) for pollster \(0\), and \(4 \to 8\) and \(3 \to 7\) for pollster \(1\). Combined walking paths and vehicular routes correspond to feasible service path that constitute the optimal daily plan. Breaks are represented by the dashed nodes \(2\) and \(3\), respectively, for pollsters \(0\) and \(1\). The optimal value is achieved within a single day using one vehicle and two pollsters, resulting in a total cost of \(\$560\).

\subsection{Formulation}

For the purpose of clarity, the integrated vehicle and pollster routing model (\modelo) is presented in the following blocks. \Cref{tb:Parameters} provides a summary of all the pertinent parameters employed in the model.

\begin{table}[htb!]
\caption{Relevant parameters and sets of the proposed model.}
\label{tb:Parameters}
\centering
\fontsize{9}{8.5}\selectfont
\setlength{\tabcolsep}{6pt}
\def\arraystretch{1.5}
\begin{tabular}{l | p{10cm}}
\toprule
	\bf Symbol & \bf Parameter
	\\
	\hline
	\(n\)		& 	Number of stores
	\\
	\( K \)	& 	Set of vehicles
	\\
	\( E \)	& 	Set of pollsters
	\\
	\( S \)	& 	Set of days in the time horizon
	\\
	\( B_{\max} \)	&	Day length
	\\
	\(P\)		&	Length of the break for pollsters
	\\
	\( Q\)		&	Vehicle capacity
	\\
	\(G \coloneqq (V,A)\) 	&  Multi-graph
	\\ 
	\(V \coloneqq C\cup \{0,2n+1\}\) 		& Node set of \(G\)
	\\
	\(C \coloneqq C_{-}\cup C_{+}\)	& Set of nodes associated with stores
	\\
	\(C_{-} \coloneqq \{1,\dots,n\}\) 	& Set of nodes representing the start of service
	\\
	\( C_{+} \coloneqq \{n+1,\dots,2n\}\) &	Set of nodes representing the end of service
	\\
	\(A \coloneqq A_S \cup A_W \cup A_V\) 		& Arc set of \(G\)
	\\
	\(A_S\)	& Service arcs connecting each store with its duplicate
	\\
	\(A_W\)	& Walking arcs
	\\
	\(A_V\)	& Vehicle arcs
	\\
	\(t_a\) for \(a\in A_S\)	& Time spent by a pollster collecting data at a store
	 \\
	 \(t_a\) for \(a\in A_W\)  & Walking travel time between stores
	 \\
	\( \tau_a\) for  \(a\in A_V\) & Vehicular travel time between stores
	 \\
	 \(\kappa_0 \) 	& Fixed operational cost per working day
	 \\
	   \(\kappa_1 \)  	& Fixed vehicle hiring cost
	   \\
	  \(\kappa_2 \) 	& Fixed pollster hiring cost
	  \\
	  \([T_0,T_1] \)	& Time window in which pollsters must start their break
	  \\
	  \(M\) & Sufficiently large number
	  \\
\bottomrule
\end{tabular}
\end{table}

\subsubsection*{Variables}

The following sets of \emph{binary} decision variables are employed in the model. First, a variable \(\xVar\) is equal to one if and only if the arc \((i, j) \in A_S \cup A_W\) is selected to be part of a walking path of pollster \(e \in E\) on day \(s \in S\). Second, a variable \(\yVar\) is only equal to one whenever arc \((i, j) \in A_V\) is included in the route of vehicle \(k \in K\) on day \(s \in S\). Third, \(\zVar\) becomes one whether pollster \(e \in E\) is transported by any vehicle over the arc \((i, j) \in A_V\) on day \(s \in S\). Fourth, \(\bVar\) (correspondingly, \(\fVar\)) is one if a pollster \(e \in E\) commences (concludes) a walking path at node \(i \in C_-\) (\(i \in C_+\)) on day \(s \in S\). Fifth, \(\wVar\) takes the value of one whenever a pollster \(e \in E\) takes a break at node \(i \in C_-\) on day \(s \in S\). Finally, a state variable \(\uVar\) is equal to one if any workload is assigned to pollsters and vehicles on day \(s \in S\) as part of the multi-day plan. Whenever any of the defining conditions of these variables are not met, then their value corresponds to zero.

In addition, for each node \(i\in C_-\), the variables \(B_i\) and \(B_{i+n}\) denote, respectively, the arrival and departure times of a pollster at nodes \(i\) and \(i+n\). These are \emph{continuous nonnegative} variables, since time is modeled as a continuous quantity bounded by \(B_{\max}\).

\subsubsection*{Objective function}

The objective function is designed to determine the minimum level of resources required to operate a feasible multi-day plan. In particular, it minimizes the total number of active operational days together with the number of vehicles and pollsters that must be hired. The first term captures the fixed daily operating cost, whereas the second and third terms represent the hiring costs of vehicles and pollsters on active days.
\begin{equation}
\label{fun:IVPRP}
	\min \,	\kappa_0 \, \sum_{s\in S} u^{s}
	+
	 \kappa_1 \, \sum_{s\in S} \sum_{k\in K} \sum_{j\in C}   y_{0,j}^{k,s}
	+ 
	\kappa_2 \, \sum_{s\in S} \sum_{e\in E} \sum_{j\in C}  z_{0,j}^{e,s}.
\end{equation}

Although uncommon in routing problems, this choice reflects the economics of the public-sector application motivating the IVPRP. In the setting considered here, the main managerial decision is not to minimize mileage or travel time per se, but to determine the smallest combination of days, vehicles, and pollsters needed to cover all stores while satisfying synchronization, shift-length, and break constraints. Accordingly, the model is intended as a \emph{resource-allocation and feasibility-planning tool} for constructing a multi-day plan, rather than as a detailed route-improvement model once such a plan has been fixed.

We acknowledge that, under a purely fixed-cost objective, two feasible solutions using the same number of days, vehicles, and pollsters may differ in terms of traveled distance or travel time. In practice, this does not contradict the purpose of the model: once the minimum-cost resource profile has been identified, the corresponding daily plans can be post-processed to improve the walking paths and vehicular routes without altering the number of active days or hired resources. This can be achieved, for example, by solving a secondary optimization problem restricted to the selected resource levels, or by applying local-search improvement procedures that minimize total travel time or distance subject to preserving feasibility and the primary cost components. Therefore, travel efficiency is treated as a \emph{secondary operational criterion}, whereas the primary optimization target is the minimum resource requirement needed to execute the survey campaign over the planning horizon.

\subsubsection*{Pollster routing}

The following set of constraints describes the routes for each pollster \(e \in E\) within each day \(s \in S\) in the time horizon.
\begin{subequations}
\begingroup
\addtolength{\jot}{0.3em}
\begin{align}
	\label{re:1a}
	\dsumin{s\in S} \sum_{e\in E} \xVar[i,i+n] &= 1		&& \forall i \in C_-,
	\\
	\label{re:1b}
	\dsumin{j \in C_+} \xVar[j,i] - \xVar[i,i+n] &= - \bVar 	 			&& \forall i \in C_-, e \in E, s \in S,
	\\
	\label{re:1c}
	\xVar[i-n,i] - \dsumin{j \in C_-} \xVar[i,j] &= \fVar[i]  	 			&& \forall i \in C_+, e \in E, s \in S,
	\\
	\label{re:1d}
	\sumin{j \in C \cup \{0\}} \zVar[j,i] 	&\leq 1 - \xVar[i,i+n] + \bVar[i]			&& \forall i \in C_-, e \in E, s \in S, 
	\\
	\label{re:1e}
	\sumin{j \in C \cup \{2n+1\}} \zVar[i,j] &\leq 1 - \xVar[i-n,i] + \fVar[i]			&& \forall i \in C_+, e \in E, s  \in S, 	
	\\
	\label{re:1f}
	\dsumin{j \in C \cup \{0\}} \zVar[j,i] \quad - \sumin{j \in C \cup \{2n+1\}} \zVar &= \bVar 		&& \forall i \in C_-,  e \in E,   s\in S,
	\\
	\label{re:1g}
	\dsumin{j\in C \cup \{0\}} \zVar[j,i]  \quad - \sumin{j \in C \cup \{2n+1\}} \zVar &=  - \fVar 		&& \forall i \in C_+,  e \in E,   s\in S,
	\\
	\label{re:1h}
	\sumin{i \in C} \zVar[0,i] &\leq u^s 			&& \forall e \in E, s\in S.
\end{align}
\endgroup
\end{subequations}
Constraints \cref{re:1a} ensure that each store \(i \in C_{-}\) is visited exactly once by a pollster on a single day. Constraints \cref{re:1b} and \cref{re:1c} delineate the edges of walking paths. Specifically, if \(\bVar = 1\) (or \(\fVar = 1\)), then a walking path must commence (or conclude) at node \(i \in C_{-}\) (or \(C_{+}\)). Furthermore, constraints \crefrange{re:1d}{re:1g} assure that if a node \(i\) is the starting or ending point of a walking path, then exactly one vehicle must pick up or drop off one pollster at that node. Finally, constraints \cref{re:1h} stipulate that each pollster departs from the depot no more than once on each active day.

\subsubsection*{Vehicle routing}

The following constraints define vehicular routes for pollster transportation.
\begin{subequations}
\begingroup
\addtolength{\jot}{0.3em}
\begin{align}
	\label{re:2a}
	\dsumin{k\in K} \sumin{j \in C \cup \{2n+1\}} \yVar &= \sumin{e\in E} \bVar 	 && \forall i \in C_-, \, \forall s \in S,
	\\
	\label{re:2b}
	\dsumin{k \in K} \sumin{j\in C\cup \{2n+1\}} \yVar &= \sumin{e \in E}  \fVar 	 && \forall i \in C_+, \, \forall s \in S,
	\\
	\label{re:2c}
	\sumin{i \in C \cup \{0\}} \yVar \quad - \sumin{i \in C \cup \{2n+1\}} \yVar[j,i] &= 0 	 && \forall j \in C,  k \in K, s \in S,
	\\
	\label{re:2d}
	\dsumin{j\in C} \yVar[0,j] - \sumin{j\in C} \yVar[j,2n+1] &= 0 	 && \forall k \in K,  s \in S,
	\\
	\label{re:2e}
	\sumin{i \in C} \yVar[0,i] 	&\leq u^s 			&& \forall k \in K, s\in S,
        \\
        \label{re:2f}
	\sumin{e\in E} \zVar 	&\leq Q \sumin{k\in K} \yVar		&&  \forall (i,j) \in A, s\in S.
\end{align}
\endgroup
\end{subequations}
Constraints \cref{re:2a} and \cref{re:2b} guarantee that vehicles arrive exclusively at pickup nodes and depart only from delivery nodes. Constraints \cref{re:2c} and \cref{re:2d} are flow conservation constraints. Additionally, \cref{re:2e} ensure that each vehicle departs from the depot at most once on each active day. Finally, \cref{re:2f} establish the connection between pollster and vehicular routes.


\subsubsection*{Time management and shift length}

The following constraints are introduced to define time bounds for the duration of a pollster shift, as well as the travel and visitation times required for each store.
\begin{subequations}
\begingroup
\addtolength{\jot}{0.3em}
\begin{align}
	&
	\label{re:3a}
	\BVar[i+n]  \geq \BVar[i] + t_{i,i+n} +  P \sumin{s\in S}\sumin{e\in E} \wVar[i]			&&  \forall i\in C_-,
	\\&
	\label{re:3b}
	\BVar[j] \geq \BVar[i] + t_{i,j}   - M \bigg(1 -  \sumin{s\in S}\sumin{e\in E} \xVar[i,j] \bigg)		&&  \forall i\in C_+, j \in C_-, 
	\\&
	\label{re:3c}
	\BVar[j] \geq \BVar[i] + \tau_{i,j}   - M \bigg(1 -  \sumin{s\in S}\sumin{k\in K} \yVar[i,j] \bigg)		&& 	\forall (i,j) \in A_V, i \neq 0, j \neq 2n+1, 
	\\&
	\label{re:3d}
	\BVar[i] \geq  \tau_{0,i} - M \bigg( 1 - \sumin{s\in S}\sumin{k\in K} \yVar[0,i] \bigg)
	&&  \forall i \in C,
	\\&
	\label{re:3e}
	B_{\max} \uVar
	\geq \BVar + \tau_{i,2n+1} - M\left( 1 - \sumin{k\in K} \yVar[i,2n+1] \right) 	&& \forall i\in C, s\in S,
\end{align}
\endgroup
\end{subequations}
Constraints \cref{re:3a} require that the departure time from a store \(i\) must be at least as large as the arrival time plus the service time and, where applicable, the time needed for a pause at that store. It is important to note that the sum in the right-most term of these constraints equals one if a pause is taken at node \(i\), and zero otherwise. Constraints \cref{re:3b} and \cref{re:3c} ensure that the arrival time at a store within a route takes into account the walking or vehicular travel time from the previously visited store, respectively. If node \(i\) immediately follows the depot in any route, constraints \cref{re:3d} establish a lower bound on the arrival time for vehicles and pollsters. Finally, constraints \cref{re:3e} link the variables \(B_i\) and \(u^s\) to bound the latest arrival time of any vehicle at the depot for each day \(s\). Here, \(M\) represents a sufficiently large number such that \(M > \max\{t, \tau\} + B_{\max}\).

\subsubsection*{Pollster breaks}

The following constraints ensure that each pollster on duty takes a break starting within the specified time frame on a given day.
\begin{subequations}
\begingroup
\addtolength{\jot}{0.3em}
\begin{align}
&
	\label{re:4a}
	T_0 \sumin{s\in S} \sumin{e\in E} \wVar \leq \BVar + t_{i,i+n} \leq   T_1  +  M \bigg( 1 - \sumin{s\in S} \sumin{e\in E} \wVar \bigg)
						&& \forall i \in C_- ,
	\\ &
	\label{re:4b} 
	\wVar \leq \xVar[i,i+n]  	&& \forall i \in C_-, e \in E, s \in S,
	\\&
	\label{re:4c} 
	\sumin{i\in C_-} \wVar = \sumin{j\in C} \zVar[0,j] 	&&\forall e \in E, s \in S.
\end{align}
\endgroup
\end{subequations}
Constraints \cref{re:4a} stipulate that if a pollster takes a break at node \(i\), it must occur within the specific time interval \([T_0, T_1]\). Constraints \cref{re:4b} guarantee that pollsters only take breaks at the stores they have visited, while constraints \cref{re:4c} enforce that each pollster on duty takes a break exactly once each day.

\subsection{Symmetry breaking inequalities}

In any feasible multi-day plan, a fixed vehicle serves a vehicular route; however, any other vehicle in \(K\) can also serve the route, yielding the same objective value. Similarly, the workload of any pair of pollsters in \(E\) can be permuted within a daily plan. Furthermore, the same principle applies to the set of days \(S\) concerning the feasible daily plans. Consequently, a set of equivalent solutions may arise from relabeling assigned vehicles and pollsters. The following set of constraints aims to prevent these symmetries. This approach can be regarded as an application of the well-known technique of symmetry elimination \citep{Kaibel-2011,Ostrowski-2009,Margot-2009}.
\begin{subequations}
\begingroup
\addtolength{\jot}{0.3em}
\begin{align}
&
	\label{re:5a}
	\sumin{i\in C} z_{0,i}^{0,0} = 1,          
	\\ &
	\label{re:5b} 
	\sumin{i\in C} \zVar[0,i] \leq \sumin{i\in C} z_{0,i}^{e-1,s} , 	&& \forall e \in E \setminus \{0\}, s\in S,
        \\ &
	\label{re:5c} 
	\sumin{i\in C} \zVar[0,i] \leq \sumin{i\in C} z_{0,i}^{e,s-1} , 	&& \forall e\in E, s \in S \setminus \{0\},
	\\&
	\label{re:5d} 
	\sumin{i\in C} y_{0,i}^{0,0} = 1,
	\\ &
	\label{re:5e} 
	\sumin{i\in C} \yVar[0,i] \leq \sumin{i\in C} y_{0,i}^{k-1,s} , 	&& \forall k \in K \setminus \{0\}, s\in S,
	\\ &
	\label{re:5f} 
	\sumin{i\in C} \yVar[0,i] \leq \sumin{i\in C} y_{0,i}^{k,s-1} , 	&& \forall k\in K, s \in S \setminus \{0\},
        \\ &
	\label{re:5g}
	\sumin{e\in E} \bVar \leq 1-\sumin{e\in E} b_{i}^{e,s-1} ,   	&& \forall i \in C_{-}, s \in S \setminus \{0\},
        \\ &
	\label{re:5h} 
	\sumin{e\in E} \fVar \leq 1- \sumin{e\in E} f_{i}^{e,s-1} , 	&& \forall i \in C_{+}, s \in S \setminus \{0\},
        \\&
	\label{re:5i}
	\dsumin{e\in E}\sumin{j\in C_+} \xVar[j,i] \leq 1 -\sumin{e\in E} b_{i}^{e,s-1}, 			&& \forall i \in C_{-}, s\in S \setminus \{0\},
	\\&
	\label{re:5j}
	\dsumin{e\in E}\sumin{j\in C_-} \xVar[i,j] \leq 1 -\sumin{e\in E} f_{i}^{e,s-1},	&& \forall i \in C_+,  s\in S \setminus \{0\},
        \\&
	\label{re:5k}
        \sum_{i \in C_-}\bVar \leq |C_-|-\sum_{r=0}^{s-1} \sum_{i\in C_-} b_{i}^{e-1,r} ,   	&& \forall e\in E \setminus \{0\}, s\in S \setminus \{0\},
        \\ &
	\label{re:5l}
        \sum_{i\in C_+}\fVar \leq |C_+|-\sum_{r=0}^{s-1} \sum_{i\in C_+} f_{i}^{e-1,r} ,   	&& \forall e\in E \setminus \{0\}, s\in S \setminus \{0\},
        \\ &
	\label{re:5m}
        u^s \leq u^{s-1} ,   	&& \forall s\in S \setminus \{0\}.
\end{align}
\endgroup
\end{subequations}
Constraints \crefrange{re:5a}{re:5c} ensure an ordered assignment of pollsters to each day \(s\in S\), starting with the pollster labeled as \(0\). Similarly, constraints \crefrange{re:5d}{re:5f} guarantee an ordered assignment for vehicles. Constraints \crefrange{re:5g}{re:5j} define daily relationships between served stores across two consecutive days. Specifically, if node \(i\) is visited on day \(s-1\), then the corresponding variables for this node on day \(s\) are constrained to zero. Constraints \crefrange{re:5k}{re:5l} impose an upper bound on the number of stores not visited in the preceding days. Finally, constraints \cref{re:5m} reinforce the activation of consecutive active days.

\subsection{Upper bounds}\label{sec:ub}

An upper bound for the objective function of {\problema} can be determined by solving the following nonlinear integer optimization problem:
\begin{subequations}
\label{NL:Upper_bound}
\begin{equation}
	\min \kappa_0 \ell_S + \kappa_1 \ell_K \ell_S + \kappa_2 \ell_E \ell_S 
	\label{eq:6a}
\end{equation}
subject to
\begin{align}
&
	\dfrac{1}{ \ell_E }  \Big( P \ell_E \ell_S + \hspace{-0.5em} \sum_{ i \in C_- } t_{i,i+n}  \Big)  
	+
	\dfrac{1}{ \ell_K }  \sum_{i \in V} \Big[ \hspace{-0.2em} \max_{(i,j) \in A_V} \tau_{i,j} + \max_{(j,i)\in A_V} \tau_{j,i} \Big]
	\leq B_{\max} \ell_S,
	\label{eq:6b}
	\\
	&\ell_K \leq \ell_E \leq Q \ell_K, \label{eq:6g}
\\
&1 \leq \ell_S  \leq |S|, 
\qquad
1 \leq \ell_K  \leq |K|, 
\qquad
1 \leq \ell_E  \leq |E|, \label{eq:6e}
\\
&\ell_S, \ell_K, \ell_E \in \mathbb{N}; \label{eq:6f}
\end{align}
\end{subequations}
where \(\ell_{S}\), \(\ell_K\), and \(\ell_E\) are decision variables representing the needed number of days, vehicles, and pollsters, respectively. The optimal solution \((\ell_S^*, \ell_K^*, \ell_E^*)\) of \crefrange{eq:6a}{eq:6e} provides an upper bound on the objective function \cref{fun:IVPRP} of the {\modelo}.  As a result, this solution determines a tighter value for the cardinality of the sets \(S\), \(K\), and \(E\) required to obtain a multi-day plan. In other words, these bounds allow us to re-define the sets \(S \coloneqq \{1,\ldots, \ell_S^*\}\), \(K \coloneqq \{1,\ldots, \ell_K^*\}\), and \(E \coloneqq \{1,\ldots, \ell_E^*\}\), which directly reduce the size of the {\modelo}. Furthermore, a solution to this problem establishes another necessary condition for the existence of feasible solutions.

Note that constraint \cref{eq:6b} divides a day into two segments. The first segment corresponds to the average daily workload per pollster, encompassing data collection activities and mandatory breaks aggregated across all days. The second segment represents the average daily vehicle routing time, calculated based on the worst-case scenario. It accounts explicitly for the maximum travel time required by vehicles, including trips undertaken solely to drop off or pick up pollsters at each node, even when vehicles travel without any pollster onboard. The sum of both terms in any feasible solution of the {\modelo} must necessarily be less than the time horizon \(B_{\max}\). Furthermore, constraint \cref{eq:6g} limits the number of pollsters to the vehicle capacity. Finally, the optimal solution can be obtained, in the worst-case scenario, using exhaustive enumeration in \(\mathcal{O}\big( |S| |K| |E| \big)\) time.

\subsection{Lower bounds}\label{sec:lb}

Similarly, it is possible to establish lower bounds based on store visiting times. To accomplish this, first, let \(\ell_e\) be the smallest integer value that specifies the minimum number of pollsters required, within a maximum of \(|S|\) days, to distribute the aggregated service time while complying with the maximum shift length and mandatory break durations. This value is determined by solving the following nonlinear optimization problem:
\begin{equation}
\label{model:lower_bounds}
	\max \bigg\{  \frac{1}{\ell_e} \sum_{ i \in C_- } t_{i,i+n} : \,
	\frac{1}{\ell_e} \sum_{ i \in C_- } t_{i,i+n} \leq B_{\max} - P	 \wedge
	\ell_e \in \{ 1, \ldots, |S||E| \}
	\bigg\}.
\end{equation}
Notice that if problem \cref{model:lower_bounds} is infeasible, then there are insufficient resources to solve the {\modelo}, thereby providing an infeasibility certificate. In other words, a solution to this problem establishes an extra necessary condition for the existence of feasible solutions.

Next, the minimum number of days required to visit all stores, denoted as \(\ell_s\), is determined as the first index \(s\) such that \(\ell_e\) belongs to the set \(\Lambda_{E,s} = \{ s, s+1, \ldots , s \cdot |E| \}\) for all \(s \in \{1, \ldots, |S|\}\). Here, \(\Lambda_{E,s}\) is the set containing the cumulative number of pollsters that can be active on a day with index \(s\). 

Finally, a bound on the number of vehicles is derived by noting that if a day is active, then there must also be an active vehicle (see \cref{re:2f}), and each vehicle can transport at most \(Q\) pollsters. Thus, the smallest number \(\ell_k\) that satisfies these constraints solves the following linear program:
\(
	\min \big\{  \ell_k : \, \ell_e \leq Q \ell_k \, \wedge \,  \ell_k \in \{ \ell_s, \ldots, \ell_s |K| \} \big\}
\).

A solution for these three problems can be obtained through exhaustive enumeration in \(\mathcal{O}\big( |S| |E| + |K| \big)\) time complexity. Moreover, while the upper bounds \((\ell_S, \ell_K, \ell_E)\) may vary based on the costs in the objective function, the lower bounds \((\ell_s, \ell_k, \ell_e)\) remain cost-independent. These lower bounds are incorporated into the {\modelo} as the following constraints:
\begin{subequations}
\begin{align}
	\sum_{s\in S} u^{s} &\geq \ell_s,
	\\
	\sum_{s\in S} \sum_{k\in K} \sum_{j\in C}   y_{0,j}^{k,s} &\geq \ell_k,
	\\ 
	\sum_{s\in S} \sum_{e\in E} \sum_{j\in C}  z_{0,j}^{e,s} &\geq \ell_e.
\end{align}
\end{subequations}
Combining these inequalities yields the lower bound \( \kappa_0 \ell_s + \kappa_1 \ell_k + \kappa_2 \ell_e \) with respect to the objective of the {\modelo}.

\section{Three-phase solution approach}\label{section_solution_approach}

Although the proposed MIP offers a comprehensive and precise representation of the IVPRP, it is not intended to be a standalone, scalable method for large real-world instances. Its rôle in this work is to serve as an initial modeling approach that captures the essence of the problem, offers an exact benchmark on small instances, and becomes a building block for developing heuristics (which also benefit from the exact bounds and symmetry-breaking inequalities). In this context, the use of heuristics is appealing for determining time- and memory-efficient solutions. This section presents an iterative method based on a graph-partitioning metaheuristic. The method comprises three sequential phases, detailed in the following subsections. Other approaches using graph-partitioning techniques have been similarly reported in the literature; see \citep{Sartori-2018,Alvarenga-2007} for further details.

As described in \cref{alg_daily_plans}, the method consists of three phases executed consecutively, following a triage principle governed by an even integer $k$.
First, the nodes are partitioned into $k$ balanced subsets using the procedure described in \cref{ssec:partitioning}, where each subset is associated with a part-time daily plan, either in the morning or afternoon, without breaks.
Second, as detailed in \cref{ssec:routing}, the graph induced by each subset is taken as input to solve a \emph{reduced} version of the main problem, referred to as ({\Rmodelo}). 
Finally, if the algorithm is successful up to this point, the resulting part-time solutions are linked to construct feasible full daily plans, see \cref{ssec:linking}. The value of $k$ is triaged against both finding a feasible $k$-partition and obtaining a feasible routing for each subset via a MIP solver; if all $k$ subproblems admit a feasible solution, the method proceeds to the linking phase, otherwise $k$ is incremented by two and the process restarts.

\begin{algorithm}[b]
\begin{algorithmic}[1]
\fontsize{9.5}{11}\selectfont
\State \textbf{Input:} an instance of the {\modelo}, maximum solver time \(T_{\max} > 0\), maximum number of iterations of the partitioning algorithm \(N_{\max} \in \mathbb{N}\), an even number \(\underline{k} \in 2\mathbb{N}\).
\State \textbf{Output:} a feasible solution of the {\modelo}.
\State \textbf{Initialization:} \(k = \underline{k}\); \(f = \texttt{false}\).
\While{\(f = \texttt{false}\)}
    \State Perform \(N_{\max}\) iterations of the partitioning procedure to find a \(k\)-way partition of \(V(G)\).
    \For{\(i \in \{1,\ldots, k\}\)}
        \State Find the best feasible solution of the {\Rmodelo} induced by subset \(i\), within the time limit \(T_{\max}\).
        \If{{\Rmodelo} induced by subset \(i\) is infeasible}
            \State \(f = \texttt{false}\)
            \State \textbf{break}
        \Else
            \State \(f = \texttt{true}\)
        \EndIf
    \EndFor
    \If{\(f\) is \(\texttt{false}\) and \(k + 2 \leq 2|S|\)}
        \State \(k = k + 2\)
    \Else
        \State \textbf{break}
    \EndIf
\EndWhile
\If{\(f\) is \(\texttt{true}\)}
    \State Perform the linking procedure.
\EndIf
\end{algorithmic}
\caption{Constructing feasible multi-day plans}
\label{alg_daily_plans}
\end{algorithm}

\subsection{Partitioning}\label{ssec:partitioning}

To divide the node set into smaller subsets, the following problem is solved: Given an undirected complete graph \(G = (C_{-}, F)\) with edge costs corresponding to the walking travel times between each pair of stores, and a weight function \(t: C_{-} \to \mathbb{R}^+\) representing the service times at stores, where \(k \geq 2\) is a fixed integer, a \(k\)-\emph{partition} of \(G\) is a collection of \(k\) subgraphs \((C_1, F(C_1)), \ldots, (C_k, F(C_k))\) of \(G\), in which \(\{ C_1, \ldots, C_k \}\) forms a partition of \(C_{-}\). Here, \(F(C)\) corresponds to the subset of edges from \(F\) with nodes in \(C\).
Moreover, let \(W_L, W_U \in \mathbb{R}^+\), with \(W_L \leq W_U\), represent the lower and upper bounds for the total service time of each subgraph, defined as the cumulative service times of all nodes within the subgraph. The problem entails finding a \(k\)-partition such that \(W_L \leq \sum_{j \in C_i} t_j \leq W_U\) for all \(i \in \{1, \ldots, k\}\), while minimizing the sum of the edge costs across all subgraphs \((C_1, F(C_1)), \ldots, (C_k, F(C_k))\).

To solve this partitioning problem, the metaheuristic proposed in \cite{Recalde-2018} is used, with the size constraints relaxed. Hence, \(W_L\) and \(W_U\) are calculated as \(\sfrac{\mu n}{k} - \sigma\) and \(\sfrac{\mu n}{k} + \sigma\), respectively, where \(\mu\) and \(\sigma\) are the sample mean and standard deviation of the service times associated with all the stores. 
The procedure is a Tabu Search algorithm based on an extension of the \(k\) equitable coloring problem for constructing feasible solutions of the balanced \(k\)-partitioning problem. The method builds a balanced partition of the vertices of a graph while minimizing total intra-group travel distance, exploring the solution space through 1-move and 2-exchange neighborhood operators with a randomized tabu tenure to avoid cycling.

\subsection{Routing pollsters and vehicles}\label{ssec:routing}

For any partition obtained in the preceding phase, a reduced version of the {\modelo} is solved to find feasible vehicular routes and walking paths for each subset. The reduced version, denoted as {\Rmodelo}, comprises of the original model with the following modifications:
\begin{enumerate}[label=\roman*.]
	\item The set of stores to be visited is given by the submultigraph induced by a subset of the partition.
	\item The total number of days simplifies to one \(( |S| = 1 )\), meaning that the reduced number of stores should be serviced by the available set of pollsters and vehicles on a single half-shift.
	\item Pollsters do not take breaks explicitly \( (w = \varnothing)\); a pause is assumed to take place at the depot once pollsters return at the end of a half-shift.
	\item The maximum working time of a half-shift discounts half of the time required for a pause; e.g., for a time window break around midday a half-shift length can be set to \( \frac{1}{2} (B_{\max} - P) \).
	\item Constraints \crefrange{re:4a}{re:4c} are removed.
	
	\item Simplified versions of the symmetry-breaking inequalities, which disregard consecutive days, are included; specifically, \crefrange{re:5a}{re:5b}, and \crefrange{re:5d}{re:5e}.
\end{enumerate}

The procedure designed in \cref{sec:ub,sec:lb} is adapted for the one-day case to obtain upper and lower bounds of the {\Rmodelo}. In particular, to improve the lower bounds, some valid inequalities are derived below.

\vspace{1\baselineskip}
\begin{remark}
	The following inequalities are valid for the {\Rmodelo}:
\begin{subequations}
\begingroup
\addtolength{\jot}{0.3em}
  \begin{align}
    \frac{1}{|E|} \sum_{i\in C_{-}} (B_{i+n}-B_i) &\leq B_{\max}, \label{valid_01}
    \\[-1.15em]
    B_{i}&\leq B_{i+n}, && \forall i\in C_{-}.  \label{valid_02}
  \end{align}
  \endgroup
\end{subequations}
\end{remark}
\noindent 
Inequality \cref{valid_01} is valid because the latest arrival time at the depot is greater than the minimum service time required to collect all data employing the maximum number of pollsters; that is, \(\sum_{i \in C_{-}} {t_{i,i+n}} \big/ {|E|} \leq B_{\max}\). Conversely, the total time spent at each store is greater than or equal to the total time required for data collection; that is, \(\sum_{i \in C_{-}} t_{i,i+n} \leq \sum_{i \in C_{-}} (B_{i+n} - B_i)\).

Although constraints \crefrange{re:3b}{re:3e} imply \cref{valid_02}, they are not necessarily valid when fractional solutions of the linear relaxation contain subcycles. Hence, as in a feasible solution, it is only possible to visit node \(i+n\) after visiting \(i\), and \(2n+1\) must be the last visited node. Therefore, \cref{valid_02} holds.

\subsection{Linking partitions}\label{ssec:linking}

Given that each partition corresponds to a part-time daily plan, the objective of this phase is to pair two of these subsets into a full-time feasible daily plan. It is important to note that, since each part-time daily plan may involve a different number of active vehicles and pollsters, it can be assumed that any vehicles or pollsters not needed for the afternoon shift can remain at the depot. Conversely, if the workforce is insufficient to fulfill the shift, additional vehicles or pollsters may join the workforce. To ensure consistency, breaks are assigned at the depot between each half-shift, as each day can be divided into two half-shifts by establishing a new time frame upper bound \(\tilde{B}_{\max} = (B_{\max} - P)/2\).

Therefore, the following two cases must be considered: (a) the required number of vehicles or pollsters differs between at least two sets, and (b) the required number of vehicles and pollsters is the same across all sets in the partition. In the former case, the cost-differing sets are linked by matching their costs in descending order in pairs. In the latter case, let \(\mathcal{K} = (V, E)\) be a complete graph, where each node in \(V = \{1, \dots, k\}\) is associated with a feasible solution of the {\Rmodelo} for any subset obtained during the partitioning step. A cost for each edge \(\{i, j\} \in E\) is defined as the sum of the arrival times at the depot for the last vehicle in partitions \(i\) and \(j\). Then, a minimum-cost matching in \(\mathcal{K}\) is found. It is noteworthy that if ties arise in case (a), they are resolved as a special instance of case (b). In either case, a feasible multi-day plan is returned.

\vspace{1\baselineskip}
\begin{remark}
	The time complexity of \cref{alg_daily_plans} can be expressed in terms of each phase of the algorithm. For each even \(k\) in the integral range between \( \underline{k}\) and \( 2|S|\) it holds that:
	\begin{enumerate}
		\item Partitioning: A feasible \(k\)-way partition can be found in at most $O(N_{\max} n^2)$.
		\item Routing: The MIP solving time for {\Rmodelo} is bounded, thus requiring at most \( k T_{\max}\) time units.
		
		\item Linking: A simple cost matching routine consists on sorting the costs of each half-shift, and then pairing the highest and smallest costs, which can be done in \( O(k \log k) \).
	\end{enumerate}
	The worst-case time complexity of the method becomes the cumulative sum of the aforementioned costs for all \(k\).
\end{remark}

\section{Computational experiments}\label{section_comp_exp}

In this section, computational experiments are presented. These consist of a set of randomly generated instances used to evaluate \cref{alg_daily_plans} against attempting to solve the Mixed-Integer Linear Program (MIP) formulation of the {\modelo} exactly. Additionally, the resolution of a real-world instance is addressed at the end of the section.

The MIP formulations of the {\modelo} and the {\Rmodelo} have been implemented in Gurobi 12 [ARM] \cite{gurobi} using the Python programming language. The tests were performed on a MacBook Pro with an M1 processor and 16GB of memory. All instances, along with the source code of the implementation, are available for download from the dedicated repository:

\vspace{0.5\baselineskip}
\begin{center}
    \noindent \href{https://github.com/andresrmt/Integrated_Routing}{\color{siamc}\texttt{https://github.com/andresrmt/Integrated\_Routing}}
\end{center}

\subsection{Simulated instances}

Random instances were generated using data from the OpenStreetMaps API via the \texttt{OSMnx} package \citep{Boeing2017}. Several Points of Interest (e.g., entertainment venues and shops) were sampled from an urban neighborhood, covering a radius of \(1\,000\) m from its center. For each pair of stores, the shortest vehicular route and shortest walking path between them were computed. The average travel time between stores was subsequently calculated using an average speed of \(30\) km/h for vehicles and \(5\) km/h for pollsters. Since vehicles adhere to street network directions, the resulting vehicle travel time matrix is asymmetric. Price surveying times, measured in minutes, were also randomly generated from a uniform distribution \(\mathcal{U}(1, 20)\).

As presented in \cref{Simulated_instances}, ten instances were simulated with varying values of \(n\), \(B_{\max}\) in minutes, \(P\) in minutes, \([T_0,T_1]\) in minutes, and \(Q\). The number of available vehicles, pollsters, and days to solve each instance was computed using the nonlinear model \crefrange{eq:6a}{eq:6f}. This was done by initializing \( |K| = 3\), \(|E| = 5\), and \(|S| = 15\), corresponding to the available real-world resources. Finally, the costs for the objective function were set to \(\kappa_0 = \$200\), \(\kappa_1 = \$100\), and \(\kappa_2 = \$40\).

\begin{table}[!htb]
\caption{Parameters associated with simulated instances.}
\label{Simulated_instances}
\centering
\fontsize{8}{8.5}\selectfont
\setlength{\tabcolsep}{6pt}
\def\arraystretch{1.75}
\begin{tabular}{l | c | cccccccccc}
\toprule
	\bf Parameter & \bf Symbol & \multicolumn{10}{c}{\bf Values}
	\\
	\hline
	\bf $\#$ Stores  &$n$ 		& 10 & 12 & 14 & 16 & 18 & 20 & 25 & 30 & 40 & 50
	\\
	{\bf $\#$ Vehicles} & $|K|$	& 3 & 2 & 2 & 2 & 2 & 3 & 2 & 2 & 2 & 3
	\\
	{\bf $\#$ Pollsters} & $|E|$	& 3 & 2 & 2 & 3 & 3 &4 & 2 & 2 & 4 & 5
	\\
	{\bf $\#$ Days} & $|S|$ 	& 1 & 2 & 2 & 2 & 2 & 1 & 2 & 2 & 2 & 2
	\\
	\bf Day length & $B_{\max}$ 	& 100 & 100 & 100 & 100 & 100 & 150 & 150 & 200 & 200 & 250
	\\
	\bf Pause length & $P$ 			& 10 & 10 & 10 & 10 & 10 & 15 & 15 & 20 & 20 & 25
	\\
	\multirow{ 2}{*}{\bf Pause time window} & \(T_0\) & 40 & 40 & 40 & 40 & 40 & 50 & 50 & 60 & 60 & 75
	\\
	 	& \(T_1\) & 60 & 60 & 60 & 60 & 60 & 100 & 100 & 140 & 140 & 175
	\\
	\bf Capacity & $Q$			& 1 & 1 & 1 & 2 & 2 & 2 & 3 & 3 & 4 & 4
	\\
\bottomrule
\end{tabular}
\setlength{\parskip}{2pt}\fontsize{9}{10}\selectfont
{The value of $|K|$ was bounded by 3, the value of $|E|$ by 5, and the value of $|S|$ by 15.
}
\end{table}

Four different strategies are employed to solve the {\problema}, each focusing on a distinct approach to the problem. The first strategy, labeled \(\mathcal{S}_1\), involves solving the MIP formulation of the {\modelo} using Gurobi in its default configuration. This B\&B-based method runs for maximum of \(3\,000\) seconds and aims to find optimal solutions whenever it is possible. The second procedure, labeled \(\mathcal{S}_2\), combines the B\&B approach with the \emph{NoRel} heuristic of Gurobi. The \emph{NoRel} heuristic is a general-purpose technique for MIPs that bypasses solving the LP relaxation and instead uses methods such as rounding and presolve reductions to identify feasible solutions. The heuristic is executed for \(120\) seconds within the solver time-limit of \(3\,000\) seconds, with the aim of generating initial solutions. The third approach is the problem-tailored heuristic described in \cref{alg_daily_plans} with \(T_{\max} = 600\) seconds, this is labeled as \(\mathcal{S}_3\). This approach is specifically designed for the structure of the {\modelo} and is particularly effective for large-scale instances due to its partitioning step.
The fourth and final strategy is the \emph{Fixing Variables Heuristic}, \(\mathcal{S}_4\). This method combines elements of the exact approach and the tailored heuristic. 
It begins after running \(\mathcal{S}_3\), which already provides a multi-day plan where breaks are taken at the depot.  To refine this solution, subsets of the binary variables associated with each half-shift (particularly, those related to the families \(x,y,z,b,f\)) are fixed at their current values and passed as inputs for the {\modelo}. Variables associated to the depot nodes \(\{0,2n+1\}\) are deliberately left unfixed, as these govern daily activation and depot-return decisions that must remain flexible across the full planning horizon. With this partial fixing in place, the model is solved within a runtime limit of \(T_{\max} = 600\) seconds to search for an integer-feasible solution. This strategy is best suited for smaller instances, as it involves solving the {\modelo}, which does not exhibit scalability. In general, feasibility of the {\modelo} is not guaranteed: one may construct tight time windows that can conflict with the break assignment inherited from \(\mathcal{S}_3\). Should this occur, resuming from Phase 1 with a new partition is sufficient to recover feasibility.
The outcomes of applying strategies \( \mathcal{S}_1, \mathcal{S}_2, \mathcal{S}_3\), and \( \mathcal{S}_4\) to the instances of \cref{Simulated_instances} are summarized in \cref{Results_simulation_Gurobi,Results_simulation_heuristics}.


\begin{table}[!hbt]
\caption{Results for the solution strategies \(\mathcal{S}_1\) and \(\mathcal{S}_2\): Gurobi B\&B and NoRel}
\label{Results_simulation_Gurobi}
\fontsize{8}{8.5}\selectfont
\setlength{\tabcolsep}{2pt}
\def\arraystretch{1.75}
\begin{tabular}{>{\centering\arraybackslash}p{2em} cccccccccccc}
\toprule
\multirow{2}{*}{$n$} & \multicolumn{3}{c}{\multirow{2}{*}{\bf Bounds}} && \multicolumn{3}{c}{\multirow{2}{*}{ \(\mathcal{S}_1\) }} && \multicolumn{3}{c}{\multirow{2}{*}{ \(\mathcal{S}_2\) }}
	\\[1.5em] 
	\cmidrule{2-4}	\cmidrule{6-8} \cmidrule{10-12} 
          &   $(\ell_k, \ell_e, \ell_s)$	& Low       & Up   && 	
                $(\hat{k}, \hat{e}, \hat{s})$ & Obj. & GAP    &&            
                $(\hat{k}, \hat{e}, \hat{s})$ & Obj.   & GAP                
	\\ \hline
10	  & (2,2,1)       & 480             & 620         && (2,2,1)   & 480        & 0.0\%   && (2,2,1)    & 480         & 0.0\%   \\
12    & (2,2,1)       & 480             & 960         && (2,2,1)   & 480        & 0.0\%   && (2,2,1)    & 480         & 0.0\%   \\
14    & (2,2,1)       & 480             & 960         && (3,3,2)   & 820        & 41.5\%  && (3,3,2)    & 820         & 41.5\%  \\
16    & (1,2,1)       & 380             & 1040        && --        & --         & --      && (2,4,2)    & 760         & 50.0\%  \\
18    & (1,2,1)       & 380             & 1040        && --        & --         & --      && (3,4,2)    & 860         & 55.8\%  \\
20    & (1,2,1)       & 380             & 660         && --        & --         & --      && (2,4,1)    & 560         & 32.1\%  \\
25    & (1,2,1)       & 380             & 960         && --        & --         & --      && (2,4,2)    & 760         & 50.0\%    \\
30    & (1,2,1)       & 380             & 960         && --        & --         & --      && (2,4,2)    & 760         & 50.0\%    \\
40    & (1,3,1)       & 420             & 1120        && --        & --         & --      && (4,8,2)    & 1120        & 62.5\%  \\
50    & (2,6,2)       & 840             & 1400        && --        & --         & --      && --         & --          & -- 
\\
\bottomrule
\end{tabular}
\setlength{\parskip}{2pt}\fontsize{9}{10}\selectfont
{Each row represents the results for an instance of size $n$ described in \cref{Simulated_instances}.
The symbol -- is reported whenever the maximum runtime was reached, and no feasible solution was found.
}
\end{table}

\Cref{Results_simulation_Gurobi} compares the behavior of the Gurobi solver under different configurations as previously outlined. The first column presents the number of nodes in the instances. The following three columns, grouped under the header \textbf{Bounds}, report the bounds obtained following the procedures described in \cref{sec:ub,sec:lb}. 
The other two multicolumns, \(\mathcal{S}_1\) and \(\mathcal{S}_2\), correspond to the results for Gurobi B\&B with and without the NoRel Heuristic.
Each consist of three columns displaying the triplet \((\hat{k}, \hat{e}, \hat{s})\), the objective value, and the gap of the solution. In this context, the triplet \((\hat{k}, \hat{e}, \hat{s})\) represents the cumulative number of vehicles used, the number of active pollsters, and the number of days in the plan, respectively. These quantities are computed as follows:
\[
	\hat{s} \coloneqq \sum_{s \in S} u^s, 		\qquad 
	\hat{k} \coloneqq \sum_{s \in S} \sum_{k \in K} \sum_{j \in C} y_{0,j}^{k,s}, 		\qquad 
	\hat{e} \coloneqq \sum_{s \in S} \sum_{e \in E} \sum_{j \in C} z_{0,j}^{e,s}.
\]
The gap is computed using the formula:
\begin{equation}\label{eq:GAP}
	\mathrm{GAP} = \frac{| \mathtt{Objective} - \mathtt{Low} |}{\mathtt{Objective}},
\end{equation}
where \( \mathtt{Objective} \) is the value of \cref{fun:IVPRP}, and \( \mathtt{Low} \) is the maximum between the lower bound obtained from \cref{sec:lb} and the lower bound obtained from the B\&B process. This convention follows that of Gurobi \cite{gurobi}, ensuring direct comparability between the gaps reported by the solver and those computed for the heuristic solutions, despite differing from the standard academic definition, which normalizes by \( \mathtt{Low} \) instead.
In terms of solution time, the only two instances for which the solver did not consume the maximum allocated time correspond to \( n \in \{10,12\}\). For these two instances, Gurobi solved the problem to optimality in less than 2 minutes.
The results showcase that the customized invocation of the \emph{NoRel} heuristic significantly contributed to the identification of feasible solutions in all instances except for one (\(n=50\)), while the Gurobi solver applied to the MIP implementation of the {\modelo}  was not able to find any starting feasible solutions on its own for \(n> 14\), even when using an HPC server.
These challenges prompted the development and application of the three-phase heuristic approach outlined in  \cref{alg_daily_plans}.


\begin{table}[!hbt]
\caption{Results of applying solution strategies \( \mathcal{S}_3 \) and \( \mathcal{S}_4 \) to the simulated instances.}
\label{Results_simulation_heuristics}
\fontsize{8}{8.5}\selectfont
\setlength{\tabcolsep}{2pt}
\def\arraystretch{1.75}
\begin{tabular}{>{\centering\arraybackslash}p{2em} cccccccccccc}
\toprule
\multirow{2}{*}{$n$} & \multicolumn{3}{c}{\multirow{2}{*}{\bf Bounds}} && \multicolumn{3}{c}{\multirow{2}{*}{ \( \mathcal{S}_3 \) }} && \multicolumn{3}{c}{\multirow{2}{*}{ \( \mathcal{S}_4 \) }}  \\[-0.5em] 
	\\
	\cmidrule{2-4}	\cmidrule{6-8} \cmidrule{10-12} 
          &   $(\ell_k, \ell_e, \ell_s)$	& Low       & Up   && 	          
                $(\hat{k}, \hat{e}, \hat{s})$ & Obj.   & GAP      &&          
                $(\hat{k}, \hat{e}, \hat{s})$ & Obj.  & GAP     
	\\ \hline
10	& (2,2,1)       & 480             & 620           && (2,2,1)    & 480         & 0.0\%    && (2,2,1)  & 480        & 0.0\%   \\
12    & (2,2,1)       & 480             & 960         && (2,2,1)    & 480         & 0.0\%    && (2,2,1)  & 480        & 0.0\%   \\
14    & (2,2,1)       & 480             & 960         && (4,4,2)    & 960         & 50.0\%   && (3,3,2)  & 820        & 41.5\%  \\
16    & (1,2,1)       & 380             & 1040        && (2,4,2)    & 760         & 50.0\%   && (2,3,1)  & 520        & 26.9\%  \\
18    & (1,2,1)       & 380             & 1040        && (2,4,2)    & 760         & 50.0\%   && (2,4,2)  & 760        & 50.0\%  \\
20    & (1,2,1)       & 380             & 660         && (2,3,1)    & 520         & 26.9\%   && (2,3,1)  & 520        & 26.9\%  \\
25    & (1,2,1)       & 380             & 960         && (2,3,2)    & 720         & 47.2\%   && (2,3,2)  & 720        & 47.2\%  \\
30    & (1,2,1)       & 380             & 960         && (2,3,2)    & 720         & 47.2\%   && (2,3,2)  & 720        & 47.2\%  \\
40    & (1,3,1)       & 420             & 1120        && (2,5,2)    & 800         & 47.5\%   && (2,5,2)  & 800        & 47.5\%  \\
50    & (2,6,2)       & 840             & 1400        && (3,7,2)    & 980         & 14.3\%   && (3,7,2)  & 980        & 14.3\% 
\\
\bottomrule
\end{tabular}
\setlength{\parskip}{2pt}\fontsize{9}{10}\selectfont
{Each row represents the results for an instance of size $n$ described in \cref{Simulated_instances}.}
\end{table}

\Cref{Results_simulation_heuristics} evaluates the performance of the two heuristic methods, namely, \( \mathcal{S}_3 \) and \( \mathcal{S}_4 \). The instance information and the meaning of each multicolumn follow the same reasoning as in \cref{Results_simulation_Gurobi}.
The results reveal several noteworthy insights. 
First, both methods found feasible solutions for all instances. 
Second, when comparing \( \mathcal{S}_2 \) against \( \mathcal{S}_3 \) and \( \mathcal{S}_4 \), the former outperforms \( \mathcal{S}_3 \) only in instance 14. For all other instances, the proposed heuristics attain better objective values when compared to \( \mathcal{S}_2 \). It is noteworthy that the \( \mathcal{S}_4 \) exhibits a marginal improvement (only two out of ) when compared to \( \mathcal{S}_3 \).

Observe that the \emph{Low} and \emph{Up} columns reported in \Cref{Results_simulation_Gurobi,Results_simulation_heuristics} can also be used to derive an upper estimate for the relative optimality gap, by replacing the objective value in \cref{eq:GAP} with the corresponding \emph{Up} bound (with values between \(23\%\) and \(63\%\)). In this way, the reported bounds provide additional information on the expected quality of the solutions.
In particular, the best objective values reported across the four solution strategies consistently yield a GAP strictly below this upper estimate, confirming that the solutions remain within a controled range relative to the bounding interval. This effect is especially visible for the instance with \(n=50\) stores, where the reported GAP is noticeably smaller than in the other nonzero cases. This is explained by the fact that, in this instance, the lower and upper bounds define a comparatively tighter interval, resulting in a stronger estimate of solution quality.
Moreover, the lower bounds coincide with the exact optimal value for the instances with \(n\in\{10,12\}\) stores, showing that the bounding procedure can be tight for certain instance configurations.

It is important to notice that the solution times of strategies \( \mathcal{S}_3 \) and \( \mathcal{S}_4 \) depend on the first two phases of \cref{alg_daily_plans}, where a node partition is found and then each subset of the partition is used to build an instance of {\Rmodelo}.
Since the number of nodes was small, the Tabu Search performed in the first phase of \cref{alg_daily_plans} took less than a second for all instances and partition sizes. For each instance, the value of \(\underline{k}\) was fixed to \(2\ell_s\).
A detailed report on the computational results of the second phase of \cref{alg_daily_plans} (routing phase) is presented in \cref{Solving-reduced-model}. Each set of stores is partitioned into two or four subsets, labeled A through D. For each subset, the table includes the number of vehicles used \(\bar{k}\), the number of active pollsters \(\bar{e}\), the cost of the half shift, and the optimality gap. 
The first two quantities are computed as follows:
\[
	\bar{k} \coloneqq \sum_{k \in K} \sum_{j \in C} y_{0,j}^{k}
	\qquad \text{and} \qquad
	\bar{e} \coloneqq \sum_{e \in E} \sum_{j \in C} z_{0,j}^{e}.
\]
Only partitions for which the {\Rmodelo} had a feasible solution for all subsets are included. A final column provides the cumulative solution times (in seconds) required to find all vehicle and pollster routes for each half-shift. This time dominates strategies \( \mathcal{S}_3 \) and \( \mathcal{S}_4 \), since the linking phase took less than a second for the former and the additional MIP for the latter was solved in about one minutes for most instances except for \(n=20\), where it took about 5 minutes.
Altogether, for \(n\geq 14\), strategies \( \mathcal{S}_3 \) and \( \mathcal{S}_4 \) showcase a competitive behavior against strategies \( \mathcal{S}_1 \) and \( \mathcal{S}_2 \).

\begin{table}[!htb]
\caption{Results of second phase of \cref{alg_daily_plans} (Routing pollsters and vehicles) applied to each subset of a partition from the instances in \cref{Simulated_instances}.
\hspace{2cm}%
}
\label{Solving-reduced-model}
\fontsize{8}{8.5}\selectfont
\setlength{\tabcolsep}{2pt}
\def\arraystretch{1.75}
\begin{tabular}{>{\centering\arraybackslash}p{2em} ccccccccccccccc cr}
\toprule
\multirow{2}{*}{$n$} & \multicolumn{3}{c}{\bf Subset A} &  & \multicolumn{3}{c}{\bf Subset B} &  & \multicolumn{3}{c}{\bf Subset C} &  & \multicolumn{3}{c}{\bf Subset D} && \multicolumn{1}{c}{\multirow{2}{*}{\bf Time} }
\\ 
\cmidrule{2-4} \cmidrule{6-8} \cmidrule{10-12} \cmidrule{14-16}  
	&   	$(\bar{k}, \bar{e})$      	& Obj.   & GAP     && 
		$(\bar{k}, \bar{e})$     	& Obj.   & GAP     &&
		$(\bar{k}, \bar{e})$         	& Obj.   & GAP     &&
		$(\bar{k}, \bar{e})$         	& Obj.   & GAP     
	\\ \hline
10	& (2,2)   & 280         & 0.0\%   &  & (2,2)   & 280         & 0.0\%	&&& &&& &&& &  13.02
	\\
12	& (2,2)   & 280         & 0.0\%   &  & (2,2)   & 280         & 0.0\%	&&& &&& &&& &  1.80
	 \\
14	& (2,2)   & 280         & 0.0\%   &  & (2,2)   & 280         & 0.0\%   &  & (2,2)   & 280         & 0.0\%   &  & (2,2)   & 280         & 0.0\% 		&& 3.16  \\
16	& (1,2)   & 180         & 0.0\%   &  & (1,2)   & 180         & 0.0\%   &  & (1,2)   & 180         & 0.0\%   &  & (1,2)   & 180         & 0.0\%  		&& 302.81  \\
18	& (1,2)   & 180         & 0.0\%   &  & (1,2)   & 180         & 0.0\%   &  & (1,2)   & 180         & 0.0\%   &  & (1,2)   & 180         & 0.0\% 		&& 111.07  \\
20	& (2,2)   & 280         & 0.0\%   &  & (2,3)   & 320         & 43.8\% 	&&& &&& &&& &  415.54 
	\\
25	& (1,2)   & 180         & 0.0\%   &  & (1,1)   & 140         & 0.0\%   &  & (1,1)   & 140         & 0.0\%   &  & (1,2)   & 180         & 22.2\% 	&& 579.07  \\
30	& (1,1)   & 140         & 0.0\%   &  & (1,1)   & 140         & 0.0\%   &  & (1,2)   & 180         & 22.2\%  &  & (1,1)   & 140         & 0.0\%		&&	448.34   \\
40	& (1,2)   & 180         & 0.0\%   &  & (1,3)   & 220         & 18.2\%  &  & (1,2)   & 180         & 0.0\%   &  & (1,2)   & 180         & 0.0\%  	&&	709.54	 \\
50	& (1,3)   & 220         & 18.2\%  &  & (2,3)   & 320         & 43.8\%  &  & (1,4)   & 260         & 30.8\%  &  & (1,3)   & 220         & 18.2\%  &&  2\,160.97
\\
\bottomrule
\end{tabular}
\end{table}


\subsection{A real-world case study}

\begin{figure}[h]
	\centering
	\includegraphics[width=\textwidth]{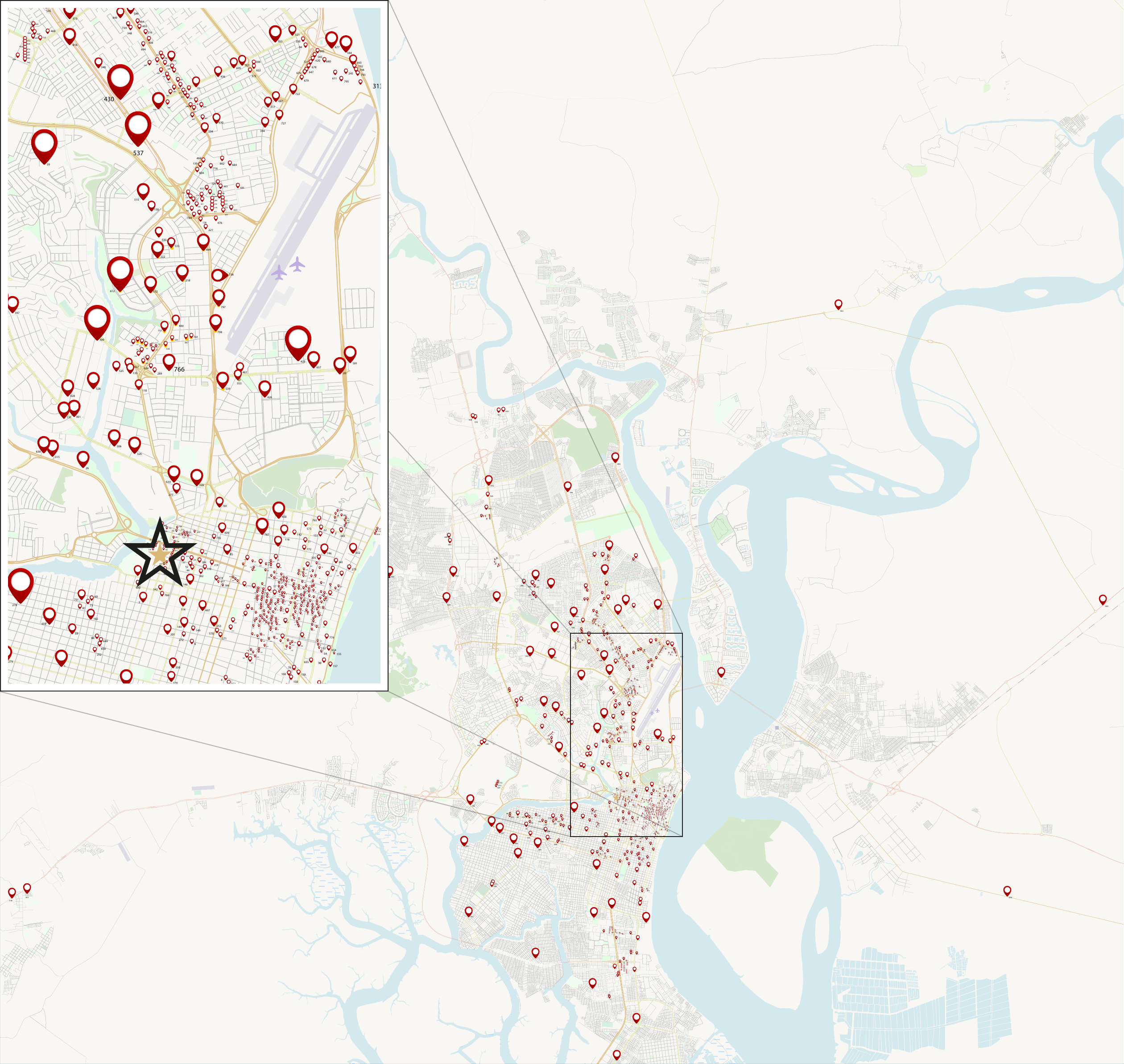}
	\caption{Map of Guayaquil showing the locations of stores to be visited.
	The baseline map was produced at a scale of 1:90\,000 (1 mm on the map represents 90 m on the ground). The inset in the upper-left corner is shown at a scale of 1:27\,000. The star highlighted in the inset marks the location of the depot (INEC headquarters). Red locator pins indicate store locations, with symbol size proportional to service time at each location.
	Note that, in addition to the fact that most stores are concentrated near the Guayas River, there is no apparent pattern in their distribution. In particular, the depot is located in the city center and is surrounded by two-high density areas (north and east) of stores.
	}
	\label{fig_mapa}
\end{figure}

The case study proposed by INEC focuses on the city of Guayaquil, the primary commercial port of Ecuador. Given its significance in economic and population development, with a population of 3 million, Guayaquil is an ideal location for data collection to monitor the CPI. Hence, INEC deploys a group of up to \(5\) pollsters and \(3\) vehicles each month to collect data from \(820\) stores within Guayaquil and its surrounding areas (as depicted in \cref{fig_mapa}). Currently, the planning process is conducted empirically, using all available resources, and data collection is completed within \(17\) days.

The operational managers at INEC aim to reduce the associated costs of the number of vehicles and pollsters within a maximum time horizon of \(15\) days. Furthermore, each day, the workforce must complete their duties within a \(7\)-hour limit, which includes a \(40\)-minute break. Accordingly,  the baseline parameters and sets are given by \(n = 820\), \(S = \{0, \ldots, 14\}\), \(E = \{0, 1, 2, 3, 4\}\), \(K = \{0, 1, 2\}\), \(B_{\max} = 420\) minutes, and \(P = 40\) minutes. Additionally, the cost parameters are set to \(\kappa_0 = \$200\), \(\kappa_1 = \$100\), and \(\kappa_2 = \$40\).

The travel time between each pair of stores was computed using the shortest path in the road and walking networks provided by a geographic information system (GIS). In the walking network, the average walking speed of pollsters is fixed at \(5\) km/h, whereas the average vehicle speed in the city is fixed at \(30\) km/h. In contrast to the simulated instances, a speed of \(40\) km/h is fixed for the outskirts of the city. Service times for each store were measured and provided by INEC. 
As outlined in \cref{sec:lb}, the computation of the lower bounds for the resulting instance of the {\modelo} indicates that any feasible multi-day plan requires at least \(3\) days, \(15\) pollsters, and \(4\) vehicles.

%
%
%
%
%
The structure of this instance makes solution strategies $\mathcal{S}_1$, $\mathcal{S}_2$, and $\mathcal{S}_4$ impractical, as the MIP formulation for the {\modelo} grows polynomially in $(n, S, E, K)$. This is unsurprizing given that the complete model contains nearly 400 million variables, rendering the resulting LP intractable due to memory constraints. Consequently, the solution strategy $\mathcal{S}_3$ was employed to obtain feasible daily plans, as its graph-partitioning phase enables it to scale effectively with instance size.

%
%

\subsubsection*{Results and solution visualization}


The following tables provide a summary of the behavior of \( \mathcal{S}_3 \) using real-world data. Notably, a lower bound on the number of partitions is given by the ratio of the total service time to the total working time available per pollster each day (\(\underline{k} = 10\)).

\begin{table}[h]
\caption{Information about the subsets of stores in the partition obtained with \cref{alg_daily_plans}.
}
\label{Statistics_instances}
\centering
\fontsize{8}{8.5}\selectfont
\setlength{\tabcolsep}{3pt}
\def\arraystretch{1.5}
\begin{tabular}{c c  c  c c c}
\toprule
{\bfseries $\#$ Subsets}  &	\multirow{2}{*}{$\mathbf{n}_{\textbf{min}}$}	&	\multirow{2}{*}{$\mathbf{n}_{\textbf{max}}$}
& {\bf Minimum}	
& {\bf Maximum}	
& {\bf $\#$ Subsets with }
\\[-0.25em]
$k$	& && {\bf service time} & {\bf service time}& {\bf feasible solutions}
\\
\hline
10 & 56     & 100    & 538.7 & 570.6 & 0   \\
12 & 36     & 97     & 446.5 & 478.1 & 0   \\
14 & 16     & 79     & 380.7 & 412.1 & 0   \\
16 & 35     & 65     & 330.9 & 362.6 & 0   \\
18 & 14     & 68     & 292.8 & 324.0  & 0   \\
20 & 11     & 60     & 261.6 & 293.2 & 1   \\
22 & 10     & 57     & 236.2 & 268.1 & 3   \\
24 & 11     & 54     & 215.2 & 247.0   & 7   \\
26 & 10     & 52     & 197.5 & 229.1 & 6   \\
28 & 7      & 43     & 134.9 & 214.1 & 16  \\
\bf 30 & \bf 7      & \bf 43     & \bf 153.7 & \bf 216.4 & \bf 30	\\
\bottomrule
\end{tabular}
\end{table}

\Cref{Statistics_instances} summarizes general information about the subsets of stores in each partition. For each \(k\) (number of subsets in a partition) and \(N_{\max} = 5\,000\), the minimum number of vertices (\(n_{\min}\)), the maximum number of vertices (\(n_{\max}\)), and the minimum and maximum service times across all subsets in the partition are reported. The final column indicates the number of subsets in the partition for which a feasible solution was found by the Gurobi solver within \(T_{\max} = 1\) hour. 
Observe that, as \(k\) increases from \(10\) to \(28\), phase two of \cref{alg_daily_plans} is unable to construct a multi-day plan. This occurs because at least one subset in the partition fails to generate a feasible vehicle-pollster route. Feasibility is achieved only when \(k=30\) (see the bold entry in \Cref{Statistics_instances}), where a feasible solution is finally obtained.

Despite the maximum number of pollsters and vehicles available at INEC being \(|E| = 5\) and \(|K| = 3\), \cref{alg_daily_plans} consistently reported that two pollsters and one vehicle are sufficient to construct every half-shift of a feasible daily plan. Consequently, \cref{30-way-part} presents the time of the last vehicular arrival time at the depot for each subset in the reported partition with \(k = 30\). 
For instance, the workload for subset \(24\) is completed in at most \(2\) hours. 
Additionally, the table marks with $\diamond$ those subsets that generated an instance of {\Rmodelo} solved to optimality within the maximum allocated time \(T_{\max}\). Moreover, Gurobi found an optimal solution for half of the \(30\) instances in less than two minutes, and only feasible solutions for the remaining instances.

\begin{table}[h]
\caption{Last vehicular arrival time at the depot for each subset of stores in the feasible solution of the routing phase \((k=30)\).}
\label{30-way-part}
\centering
\fontsize{8}{8.5}\selectfont
\setlength{\tabcolsep}{7.5pt}
\def\arraystretch{1.5}
\begin{tabular}{lc  lc  lc c  c  c  c  c  c  c c  c c}
\toprule
\multicolumn{2}{c}{\bfseries Subsets 0 --- 9}  & \multicolumn{2}{c}{\bfseries Subsets 10 --- 19} & \multicolumn{2}{c}{\bfseries Subsets 20 --- 29}\\
\cmidrule(lr){1-2} \cmidrule(lr){3-4} \cmidrule(lr){5-6} 
 &  Arrival time &  
 &  Arrival time &  
 &  Arrival time    \\
\cmidrule(lr){2-2} \cmidrule(lr){4-4} \cmidrule(lr){6-6}
0         		& 1 h 36 min & 10 $\diamond$       	& 2 h 28 min & 20        & 3 h 03 min  \\
1         		& 1 h 55 min & 11 $\diamond$       	& 2 h 14 min & 21 $\diamond$       & 2 h 01 min  \\
2 $\diamond$  	& 2 h 02 min  & 12       			& 1 h 39 min & 22        & 1 h 43 min \\
3         		& 1 h 55 min & 13 $\diamond$       	& 2 h 19 min & 23        & 1 h 45 min \\
4 $\diamond$  	& 2 h 08 min  & 14       			& 1 h 55 min & 24 $\diamond$       & 2 h 00 min  \\
5 $\diamond$  	& 1 h 57 min & 15 $\diamond$       	& 3 h 30 min & 25 $\diamond$       & 2 h 08 min  \\
6 $\diamond$  	& 2 h 18 min & 16 $\diamond$       	& 2 h 42 min & 26 $\diamond$       & 2 h 18 min \\
7         		& 1 h 50 min & 17        			& 2 h 14 min & 27        & 2 h 18 min \\
8         		& 2 h 49 min & 18 $\diamond$       	& 2 h 16 min & 28        & 1 h 47 min \\
9 $\diamond$        & 3 h 03 min  & 19       			& 2 h 04 min  & 29       & 1 h 54 min   \\
\bottomrule
\end{tabular}
\setlength{\parskip}{2pt}\fontsize{9}{10}\selectfont
{Subsets marked with $\diamond$ correspond to instances of the {\Rmodelo} solved to optimality within the maximum allocated time \(T_{\max}\).
}
\end{table}

The outcomes of the linking partition step are presented in \cref{linking_table}. Each day \(s \in \{0, \ldots, 14\}\) is assigned to a pair of subsets, corresponding to a morning and afternoon half-shift necessary to complete the workload. For example, day \(11\) synchronizes subsets \(27\) and \(1\), resulting in a total daily working time of \(4\) hours and \(53\) minutes, including a \(40\)-minute break between shifts. The first subset of stores is surveyed in the morning, while the second is visited in the afternoon.

As shown in Table \ref{30-way-part}, each daily plan is executed by \(2\) pollsters and \(1\) vehicle, successfully completing the workload within a time horizon of \(15\) days. This outcome differs from the empirical solution provided by INEC, which required \(5\) pollsters and \(3\) vehicles within a time horizon of \(17\) days. 
Given that the daily income of each pollster is \(\$40\) and the daily cost of hiring and using a vehicle is \(\$100\), the empirical solution implemented by INEC, obtained without the aid of any optimization tool, approximately doubled the cost of the multi-day plan obtained with \cref{alg_daily_plans}.

\begin{table}[h]
\centering
\caption{Multi-day plan provided by \cref{alg_daily_plans}.
}
\label{linking_table}
\fontsize{9}{8}\selectfont
\setlength{\tabcolsep}{2pt}
\def\arraystretch{1.5}
\begin{tabular}{>{\centering\arraybackslash}p{2em}    ccc>{\centering\arraybackslash}p{6em}        c       >{\centering\arraybackslash}p{2em}  ccc>{\centering\arraybackslash}p{6em}}
\toprule
\multicolumn{5}{c}{\bfseries Daily plans 0 --- 7}  &&	 \multicolumn{5}{c}{\bfseries Daily plans 8 --- 14} 
\\
\cmidrule(lr){1-5} \cmidrule(lr){7-11}
\multirow{2}{*}{\(s\)}		 & \multicolumn{3}{c}{Synced}  & Working 	&&	 \multirow{2}{*}{\(s\)}		 & \multicolumn{3}{c}{Synced}  & Working
 \\
 &	\multicolumn{3}{c}{subsets}	&	time		&&	&	\multicolumn{3}{c}{subsets}	&	time
 \\
\cmidrule(lr){1-5} \cmidrule(lr){7-11}
	0 & 15 & $\to$ & 0 & 5 h 45 min	&& 		8 & 20 & $\to$ & 12 & 5 h 22 min		\\
	1 & 9 & $\to$ & 22 & 5 h 26 min	&& 		9 & 8 & $\to$ & 23 & 5 h 14 min		\\
	2 & 16 & $\to$ & 28 & 5 h 10 min	&& 	10 & 10 & $\to$ & 7 & 4 h 59 min		\\
	3 & 13 & $\to$ & 29 & 4 h 53 min	&& 	11 & 27 & $\to$ & 1 & 4 h 53 min		\\
	4 & 6 & $\to$ & 14 & 4 h 54 min		&& 	12 & 26 & $\to$ & 3 & 4 h 53 min		\\
	5 & 18 & $\to$ & 5 & 4 h 53 min		&& 	13 & 11 & $\to$ & 24 & 4 h 55 min		\\
	6 & 17 & $\to$ & 21 & 4 h 56 min		&& 	14 & 4 & $\to$ & 2 & 4 h 51 min		\\
	7 & 25 & $\to$ & 19 & 4 h 52 min		&& 	& \multicolumn{3}{c}{\bf Average:} & 5 h 03 min     
	\\
\bottomrule
\end{tabular}
%
\setlength{\parskip}{2pt}\fontsize{9}{10}\selectfont
{The index \(s\) corresponds to the assigned day of the feasible daily plan constructed with strategy \(\mathcal{S}_3\).
}
\end{table}

\begin{figure}[h]
	\centering
	\includegraphics[width=\textwidth]{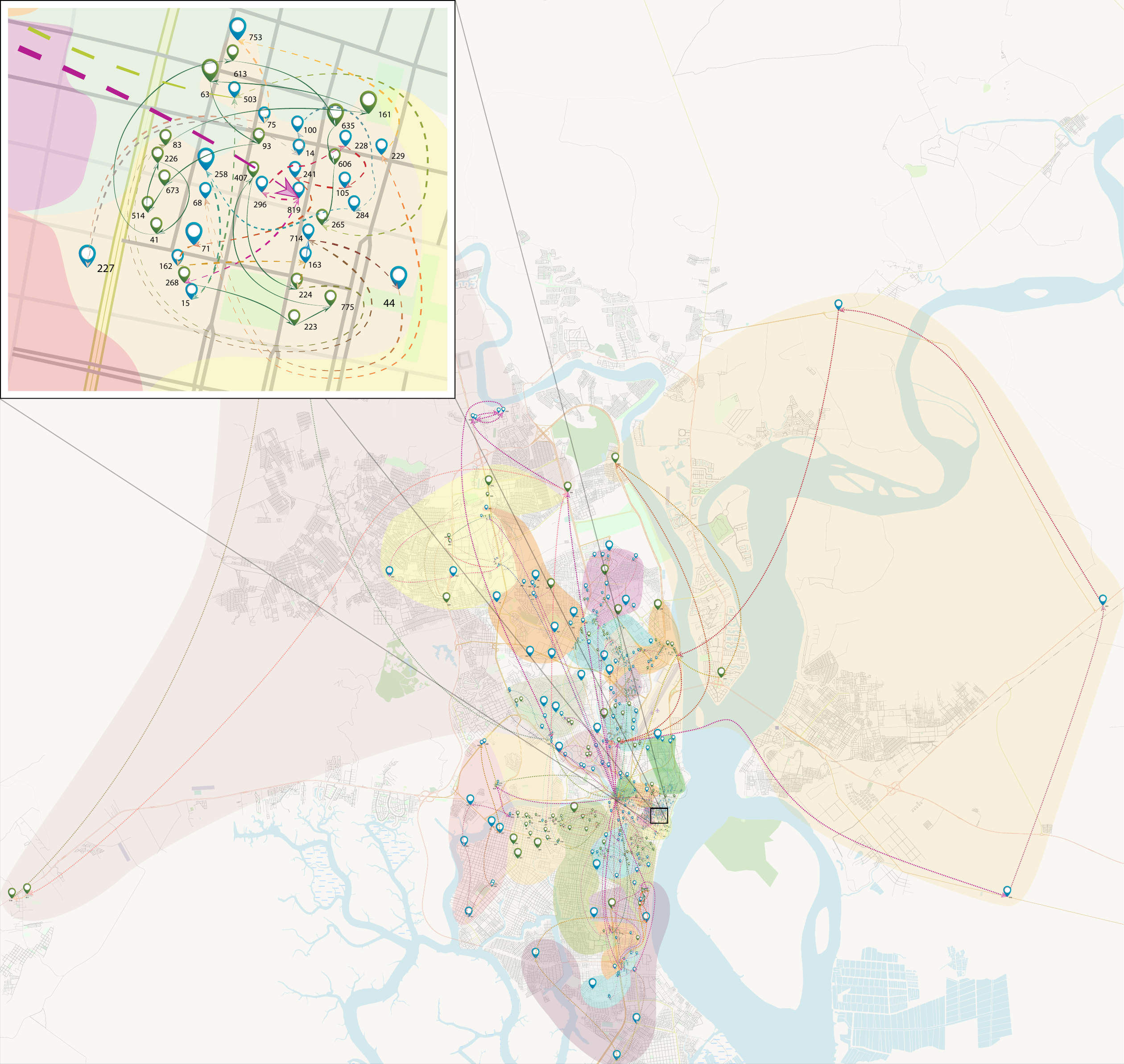}
	\caption{Map of Guayaquil illustrating the optimized partitioning of store locations into operational subsets together with the corresponding vehicle routing  and walking path structure. The base map was produced at a scale of 1:90\,000 (1 mm on the map represents 90 m on the ground). The inset shown in the upper-left corner presents a detailed view of subset~12 at a scale of 1:3\,000.
Stores are grouped into subsets, each delineated by a semi-transparent colored region representing the spatial partition assigned to a single operational tour. Teal locator pins denote vehicle-only service locations, i.e., nodes at which the pollster remains with the vehicle and no walking path is initiated. Green locator pins indicate stores visited along walking paths.
Dashed arrows represent vehicle movements between successive service locations; variations in color shade distinguish consecutive routes within each subset. Continuous green arrows indicate walking connections linking stores visited on foot. The illustrated connections represent logical routing decisions derived from the optimization model rather than the exact physical trajectories followed on the street network; actual travel between nodes corresponds to the shortest directed path under a Manhattan-distance metric.}
	\label{fig_solucion_zoom}
\end{figure}

Finally, the resulting real multi-day plan is illustrated in \cref{fig_solucion_zoom}. In the upper left corner, the afternoon shift of the eighth working day (subset 12) is enlarged. The stores visited by the first and second pollsters are represented by circles and squares, respectively, with each store labeled by a number between \(1\) and \(820\). Continuous and dashed lines indicate the routes: continuous lines represent walking paths, while dashed lines represent vehicle routes.


\section{Conclusions}\label{section_conclusions}

The Integrated Vehicle and Pollster Routing Problem (\problema), focused on the transportation of pollsters with walking paths and vehicular routes, is introduced in this paper. 
%
A mixed-integer formulation was proposed for this problem and assessed using simulated instances and real-world data. Due to the hardness of the problem, a three-phase heuristic method was devised: initially, the set of stores to be visited is partitioned, followed by the routing of pollsters and vehicles for each subset, and finally, pairs of subsets are linked as half-shifts of a multi-day plan.

The algorithm was tested with a real-world instance in Guayaquil, which is the second most important city in Ecuador.  The National Statistics Bureau of Ecuador (INEC) provided data, which included details concerning 820 stores to be visited monthly. The heuristically-found solution outperformed the one designed by INEC, approximately halving the operational costs. This was achieved by using half of the available workforce and fleet, thereby optimizing resource allocation. Moreover, the resulting multi-day plan takes two days less than the former empirical solution.

Future research could concentrate on developing more advanced heuristic and metaheuristic approaches to tackle larger instances of the {\problema}. This could involve hybrid methods that integrate exact and heuristic techniques to balance solution quality with computational efficiency. A particular choice is the incorporation of problem-specific valid inequalities, including subtour elimination inequalities, particularly within reduced formulations. In addition, alternative exact solution paradigms deserve further study. In particular, column generation, possibly embedded within a branch-and-price framework, could provide a more scalable exact methodology by using route-based rather than arc-based variables. Furthermore, there is substantial potential to explore extensions of the {\problema} that incorporate more complex constraints, such as dynamic routing adjustments based on real-time data or multi-period planning that accounts for demand variability over time. While this study focuses on a specific case in Ecuador, future research could examine how the proposed framework might be adapted to similar applications in other fields, including logistics, urban planning, and disaster response.

\section*{Declarations}

\subsubsection*{Ethics approval and consent to participate}
Not applicable.

\subsubsection*{Consent for publication}
Not applicable.

\subsubsection*{Funding} 
This project was funded by Escuela Politécnica Nacional Research Project PIJ--15-12 with support of the Laboratory of Scientific Computing of the Research Center for Mathematical Modeling -- ModeMat (\url{hpcmodemat.epn.edu.ec}). 
Andrés Miniguano Trujillo acknowledges the support of MAC-MIGS CDT Scholarship under EPSRC grant \texttt{EP/S023291/1}.

\subsubsection*{Availability of data and materials}
All instances, along with the source code of the implementation, are available for download from the dedicated repository:

\vspace{0.5\baselineskip}
\begin{center}
    \noindent \href{https://github.com/andresrmt/Integrated_Routing}{\color{siamc}\texttt{https://github.com/andresrmt/Integrated\_Routing}}
\end{center}

\subsubsection*{Competing interests}
The authors declare that they have no competing interests

\subsubsection*{Author contributions}
All authors contributed equally to this article. Under the CRediT authorship taxonomy, this includes:
Writing -- original draft, Visualization, Validation, Software, Methodology, Formal analysis, Conceptualization.

\subsubsection*{Acknowledgements} The authors are grateful with the managers of INEC, specially with Alexandra Enríquez for providing the data used to build the real-world instance, and for her valuable feedback on the results. Additionally, we acknowledge the collaboration of Pablo Zuleta during the first stage of this research.

%

\bibliography{Main.bbl}

\end{document}